\documentclass[10pt]{article}
\usepackage{graphicx}
\usepackage{amssymb}
\usepackage{epstopdf}
\usepackage{amsmath}
\usepackage{amsfonts}

\def\phi{{\varphi}}
\def\H{{\mathcal H}}
\def\P{{\mathcal P}}

\DeclareSymbolFont{AMSb}{U}{msb}{m}{n}
\DeclareMathSymbol{\N}{\mathbin}{AMSb}{"4E}
\DeclareMathSymbol{\Z}{\mathbin}{AMSb}{"5A}
\DeclareMathSymbol{\R}{\mathbin}{AMSb}{"52}
\DeclareMathSymbol{\Q}{\mathbin}{AMSb}{"51}
\DeclareMathSymbol{\I}{\mathbin}{AMSb}{"49}
\DeclareMathSymbol{\C}{\mathbin}{AMSb}{"43}

\def\be{\begin{equation}}
\def\ee{\end{equation}}
\def\ber{\begin{eqnarray}}
\def\eer{\end{eqnarray}}

\def\sv{{\bf s}}
\def\fv{{\bf f}}

\def\beq{\begin{equation}}
\def\eeq{\end{equation}}

\begin{document}

\addtolength{\textheight}{0 cm} \addtolength{\hoffset}{0 cm}
\addtolength{\textwidth}{0 cm} \addtolength{\voffset}{0 cm}

\newenvironment{acknowledgement}{\noindent\textbf{Acknowledgement.}\em}{}

\setcounter{secnumdepth}{5}
 \newtheorem{proposition}{Proposition}[section]
\newtheorem{theorem}{Theorem}[section]
\newtheorem{lemma}[theorem]{Lemma}
\newtheorem{coro}[theorem]{Corollary}
\newtheorem{remark}[theorem]{Remark}
\newtheorem{claim}[theorem]{Claim}
\newtheorem{conj}[theorem]{Conjecture}
\newtheorem{definition}[theorem]{Definition}
\newtheorem{application}{Application}

\newtheorem{corollary}[theorem]{Corollary}

\title{Symmetric Monge-Kantorovich problems and polar decompositions of vector fields}
\author{Nassif  Ghoussoub\thanks{Partially supported by a grant
from the Natural Sciences and Engineering Research Council of Canada.}\\
{\it\small Department of Mathematics}\\
{\it\small  University of British Columbia}\\
{\it\small Vancouver BC Canada V6T 1Z2}\\
{\it\small nassif@math.ubc.ca}\vspace{1mm}\and
Abbas Moameni\thanks{Research  supported by a grant
from the Natural Sciences and Engineering Research Council of Canada.}\\
\hspace{2mm}\\
{\it\small Department of Mathematics and Computer Science}\\
{\it\small University of Lethbridge}\\
{\it\small Lethbridge, AB  Canada T1K 3M4 }\\
{\it\small  abbas.momeni@uleth.ca}\\
\date{February 10, 2013, Revised September 5, 2013}
}
\maketitle

\begin{abstract} For any given integer $N\geq 2$, we show that every bounded measurable vector field from a bounded domain $\Omega$ into $\R^d$ is $N$-cyclically monotone up to a measure preserving $N$-involution. The proof involves the solution of a multidimensional symmetric Monge-Kantorovich problem, which we first study in the case of a general cost function on a   product domain $\Omega^N$. The polar decomposition described above corresponds to a special cost function derived from the vector field in question (actually $N-1$ of them). In this case, we show that the supremum over all  probability measures on $\Omega^N$ which are invariant under cyclic permutations and with a given first marginal $\mu$, is attained on a probability measure that is supported on the graph of a function of the form $x\to (x, Sx, S^2x,..., S^{N-1}x)$, where $S$ is a $\mu$-measure preserving transformation on $\Omega$ such that $S^N=I$ a.e. The proof exploits a remarkable duality between such involutions and those Hamiltonians that are $N$-cyclically antisymmetric.

\end{abstract}

\section{Introduction}  Given Borel probability measures  $\mu_1, \mu_2, ..., \mu_N$  on a domain $\Omega$ of $\R^d$ and a bounded above upper semi-continuous cost function $c: \Omega^N\to \R\cup \{-\infty\}$,  the multi-marginal version of the  Monge-Kantorovich problem
consists of maximizing
\[
\int_{\Omega^N}c(x_1,x_2,..., x_N) d\pi
\]
among all probability measures $\pi$ on $\Omega^N$ whose i-th marginal  is equal to $\mu_i$  for each $i=1,..., N$. We shall use the notation
\begin{equation}\label{gm1}
\hbox{${\rm MK}(c\, ; \mu_1,..., \mu_N)=\sup\{\int_{\Omega^N}c(x_1,x_2,..., x_N) d\pi;\, \pi \in {\cal P}(\Omega^N)$, \, ${\rm proj}_i\pi=\mu_i$ for $i=1,...,N\}$. }
\end{equation}
 In this paper we are concerned with the following symmetric counterpart of the Monge-Kantorovich problem:
 \begin{equation}\label{gm22}
{\rm MK}_{\rm sym}(c,\mu)=\sup\left\{\int_{\Omega^N}c(x_1,x_2,..., x_N) d\pi;\, \pi\in \P_{sym}(\Omega^N,\mu)\right\}
\end{equation}
 where $\P_{sym}(\Omega^N, \mu)$ denotes the set of all probability measures on $\Omega^N$,  which are invariant under the cyclical permutation
 \[ \sigma (x_1, x_2,..., x_N)=(x_2, x_3,..., x_N, x_1).
 \]
 and whose marginals are --necessarily-- equal to the same probability measure $\mu$ on $\Omega$.
 In other words, $\pi \in \P_{sym}(\Omega^N, \mu)$ if
 \begin{equation}
\hbox{$ \int_{\Omega^N}f(x_1,x_2,..., x_N) d\pi=\int_{\Omega^N}f(\sigma (x_1,x_2,..., x_N)) d\pi$ for every $f\in C(\Omega^N)$},
 \end{equation}
 and for every $i=1,..., N$,
\begin{equation}
 \hbox{$ \int_{\Omega^N}f(x_i) d\pi=\int_{\Omega}f(x_i)d\mu$ for every $f\in C(\Omega)$}.
 \end{equation}
 Standard results show that there exists $ \pi_0\in \P_{sym}(\Omega^N,\, \mu)$   where the supremum above is attained.  In this paper, we are interested in an important class of cost functions $c$, where the optimal measure $\pi_0$ is necessarily supported on the graph of a function of the form $x\to (x, Sx, S^2x,..., S^{N-1}x)$, where $S$ is a $\mu$-measure preserving transformation on $\Omega$ such that $S^N=I$ a.e.

 If $c$ is finite, then one can extend the original approach of Kantorovich to the multi-marginal and cyclically symmetric case to show that (\ref{gm22}) is dual to the following
 minimization problem
\begin{equation}\label{gm3}
{\rm DK}^1_{\rm sym}(c,\mu):=\inf \biggl\{N\int_\Omega u(x)\, d\mu; \, u:\Omega \to \R\,\, {\rm and}\,\,  \sum_{j=1}^Nu(x_j)\geq \frac{1}{N}\sum_{i=0}^{N-1}c(\sigma^i (x_1,\dots,x_N))\biggr\}.
\end{equation}
In this paper, we introduce a new dual problem involving the class ${\mathcal H}_N(\Omega)$ of all {\em $N$-cyclically antisymmetric Hamiltonians} on $\Omega^N$, that is
\begin{equation}
\hbox{${\mathcal H}_N(\Omega)=\{H\in C(\Omega^N; \R);\, \sum_{i=0}^{N-1}H\left( \sigma^{i}\left( {\bf x}\right) \right) =0\, $ for all ${\bf x}\in \Omega^N\}$.}
\end{equation}
For $H\in \H_N(\Omega)$, let $\ell_H^{^c}$ be the ``$c$-Legendre transform" of $H$ with respect to the last $(N-1)$ variables, i.e.,
\begin{equation*}
\ell_H^{^c}(x)=\sup\left\{c(x, x_2,..., x_N)-H\left( x,x_2,..., x_N\right); (x_{2},...,x_{N})\in \Omega^{N-1}\right\},
\end{equation*}
and consider the problem
\begin{eqnarray}\label{dk}
{\rm {\rm DK}}^2_{\rm sym}(c, \mu):= \inf\left\{\int_\Omega \ell_H^c(x)d\mu (x);\, H\in \H_N(\Omega)\right\}.
\end{eqnarray}
We start by proving in section 1 the following result.
\begin{theorem} \label{symmetric.dual} Let $c$ be a cost function that is continuous and bounded above, then
\begin{equation}\label{equa.1}
{\rm MK}_{\rm sym}(c,\mu)={\rm DK}^1_{\rm sym}(c,\mu)={\rm DK}^2_{\rm sym}(c,\mu).
\end{equation}
Moreover, the three extrema are attained.
\end{theorem}
Of great interest is to determine for which cost functions $c$, problem ${\rm MK}_{\rm sym}(c,\mu) $ is attained at an extremal  probability measure $\pi_0\in \P_{sym}(\Omega^N,\, \mu)$
that is  those supported on the graph of the form $x\to (x, Sx, S^2x,..., S^{N-1}x)$, where $S$ is a $\mu$-measure preserving transformation on $\Omega$ such that $S^N=I$ a.e. Indeed, it is clear that ${\rm MK}_{\rm sym}(c,\mu)\geq {\rm MK}_{\rm cyc}(c,\mu)$, where
\begin{equation}
{\rm MK}_{\rm cyc}(c,\mu):=\hbox{$\sup\large\{\int_{\Omega^N}c(x, Sx,..., S^{N-1}x) d\mu; S$ is $\mu$-measure preserving on $\Omega$ and $S^N=I$ a.e.$\large\}$}
\end{equation}
Recently, and after the first version of this paper appeared on arxive, Colombo and Di Marino established the following natural result.
\begin{theorem} (Colombo-Di Marino) Let $c$ be a cost function that is continuous and bounded above. if $\mu$ has no atoms, then
\begin{equation}
{\rm MK}_{\rm sym}(c,\mu)={\rm MK}_{\rm cyc}(c,\mu).
\end{equation}
\end{theorem}
In section 2, we  give a sufficient condition on the cost function and on the optimal anti-symmetric Hamiltonian $H$ that will insure that  ${\rm MK}_{\rm cyc}(c,\mu)$ is attained.

\begin{theorem} \label{extremal} Let $c$ be a cost function that is continuous and bounded above.
Assume that
${\rm DK}^2_{\rm sym}(c,\mu)$ 
 is attained at some $H_\infty \in \H_N(\Omega)$ in such a way  that for $\mu$-almost $x\in \Omega$ the map
 \[
 (x_2, x_3,..., x_N) \to c(x, x_2,..., x_N)-H_\infty(x,x_2,..., x_N)
 \]
  attains its maximum uniquely. Then ${\rm MK}_{\rm cyc}(c,\mu)$ is attained.
 \end{theorem}
Assume now that the cost function itself $c:\Omega^N\to \R$, is itself $\sigma$-symmetric, that is
\[
\hbox{$c(x_1, x_2,..., x_N)=c(x_2, x_3,...,x_N, x_1)$ \quad on $\Omega^N$.}
\]
The symmetric Monge-Kantorovich problem
 is then clearly equivalent to the classical one when all marginals are the same and equal to $\mu$, that is
\[
\hbox{${\rm MK}_{\rm sym}(c, \mu)={\rm MK}(c;\, \mu,...,\mu)=\sup\{\int_{\Omega^N}c(x_1, x_2,..., x_N) d\pi;\, \pi\in \P(\Omega^N)\, \, \&\, \, {\rm proj}_i\pi=\mu, i=1,..., N\}.$}
\]
Examples of such cost functions are:\\

\noindent {\bf 1. The quadratic cost} $c({\bf x})=-\sum_{i=1}^N\sum_{j=i+1}^N|x_i-x_j|^2$ considered in full generality by Gangbo-Sweich \cite{GS}. The symmetric version of the problem, that is when all the marginals are identical, can only admit the trivial solution. In other words the infimum is uniquely attained at the image of $\mu$ by the map $x\to (x,x,...,x)$, i.e., the involution $S$ is nothing but the identity. \\

\noindent {\bf 2. The Plakhov cost} function $c(x,y)=-1-\cos (x-y)$, which was studied in detail in \cite{Pla}. This is an example where ${\rm MK}_{\rm sym}(c,\mu)$ does not have a Monge solution, i.e., it  is not attained at a measure $\pi_S$ that is the image of $\mu$ by a map $x\to (x, Sx, S^2x,..., S^{N-1}x)$, where $S$ is an $N$-involution.\\

\noindent {\bf 3. The Coulomb cost}  $c({\bf x})=-\sum_{i=1}^N\sum_{j=i+1}^N\frac{1}{|x_i-x_j|}$ is a most interesting example since it appears in electronic structure theory. Indeed, recent insight into exchange-correlation in density functional theory led many authors (such as Buttazzo-De Pascale and Gori-Giorgi \cite{BDG} and Cotar-Friesecke-Kl\"uppelberg \cite{CFK}  to reformulate the electron-electron interaction energy functional with respect to a density $\rho({\bf x})$ as an $N$-dimensional mass transport, where the cost functional is the Newtonian potential
\begin{equation}
\hat V_{ee}=\sum_{i=1}^N\sum_{j=i+1}^N|{\bf x_i}-{\bf x_j}|^{-1},
\end{equation}
which correspond to $N$ interacting electrons. Assuming that the admissible configurations of $N$ electrons in $d$-dimensions
have the form $(\fv_1(\sv),....,\fv_N(\sv))$
where $\sv$ is a $d$-dimensional vector that determines the position of, say, electron ``1",
 and $\fv_i(\sv)$ ($i=1,...,N$, $\fv_1(\sv)=\sv$) are the {\it co-motion functions},  which determine the position of the $i$-th electron in terms of $\sv$, and if the variable $\sv$ itself is distributed according to the normalized density $\rho(\sv)/N$, then the energy functional  for ``strictly correlated electrons''  $V_{ee}^{\rm SCE}[\rho]$ corresponding to the density $\rho$ is given by the infimum of
\beq \label{E1}
\int_{\R^d} d\sv\, \frac{\rho(\sv)}{N} \, \sum_{i=1}^N\sum_{j=i+1}^N\frac{1}{|\fv_i(\sv)-\fv_j(\sv)|},
\eeq
among all co-motion functions $(\fv_i)_i$ that preserve the density $\rho$, so as to ensure the indistinguishability of the $N$ electrons. Formally, such functions must satisfy the equations
 \begin{equation}
\rho(\fv_i(\sv))d\fv_i(\sv)=\rho(\sv)d\sv, \, i=1, ..., N.
\label{eq_df}
\end{equation}
A relaxation of this formulation is to consider $V_{ee}^{\rm SCE}[\rho]$ as the infimum of
  \begin{equation} \label{E2}
\int_{\R^{Nd}}\sum_{i=1}^N\sum_{j=i+1}^N\frac{1}{|x_i-x_j|}d\pi (x_1,..., x_N),
\end{equation}
over all probability densities on $\R^{Nd}$ (``wave functions") whose all marginals are equal to $\frac{\rho}{N}$. The problem posed above,  that is searching for the minimum possible interaction energy in a given density, is a typical Monge-Kantorovich problem involving symmetry. The main open problem here is whether there is indeed an optimal co-motion functions  $(\fv_i(\sv)_{i=1}^N$ that minimizes both expressions in (\ref{E1}) and (\ref{E2}).
If this is the case, then our Lemma \ref{magic} below shows that the optimal one must be of the form
\[
\hbox{$\fv_1(\sv)=\sv$, $\fv_i(\sv)=f^i(\sv)$ for $i=2,..., N$, while $\fv^{N+1}(\sv)=\sv$ for some $f$ satisfying (\ref{eq_df}).}
\]
  In other words, for each density $\rho$, there exists such an $f$ satisfying (\ref{eq_df}) such that
  \[
  V_{ee}^{\rm SCE}[\rho]=\int_{\R^d} d\sv\, \frac{\rho(\sv)}{N} \, \sum_{i=1}^N\sum_{j=i+1}^N\frac{1}{|\fv^i(\sv)-\fv^j(\sv)|}.
  \]
 We note that this has been verified in the case $N=2$  \cite{BDG}. The case where there is a higher number of electrons is much more delicate \cite{CFP}.

In this paper, we shall resolve this problem for cost functions on $\Omega^N$ of the form
\begin{equation}\label{cost}
c(x_1, x_2,..., x_{N})=\langle u_1(x_1), x_2\rangle +....\langle u_{N-1}(x_1), x_N\rangle,
\end{equation}
where $u_1,..., u_{N-1}$ are given vector fields from $\Omega$ to $\R^d$. Note that this cost is not cyclically symmetric, yet we shall optimize it on the class of symmetric probabilities in order to establish an interesting representation for general vector fields in term of cyclically monotone operators. This makes use of a remarkable duality between three fundamental concepts in functional analysis: monotonicity, cyclical symmetry and involutions.
Indeed, let ${\mathcal S}(\Omega)$ denote the set of measure preserving transformations on $\Omega$, which can be considered as a closed subset of the sphere of $L^2(\Omega, \R^d)$ and let
\[
\hbox{${\mathcal S}_N(\Omega)=\huge\{S\in {\mathcal S}(\Omega), S^N=I$ $\mu$ a.e\}.}
\]
 The set ${\mathcal S}_N(\Omega)$ has been shown recently \cite{GG} to be polar to the class of {\em $N$-cyclically monotone vector fields}, which are those $u:\Omega \to \R^d$ that satisfy
 for
every cycle $x_{1},...,x_{N},x_{N+1}=x_{1}$ of points in $\Omega$, the inequality
\begin{equation}
\sum_{i=1}^{N}\left\langle u(x_{i}) ,x_{i}-x_{i+1}\right\rangle
\geq 0.
\end{equation}
More generally, Galichon-Ghoussoub \cite{GG} introduced the following extension to the case of more than one vector field.
\begin{definition} A family of vector fields  $u_1, u_2,..., u_{N-1}$ from $\Omega \to \R^d$ is said to be  {\it jointly  $N$-monotone} if for every cycle $x_1,..., x_{2N-1}$ of points in $\Omega$ such that $x_{N+l}=x_l$ for $1 \leq l \leq N-1$, we have
 \begin{equation}\label{uni0} \sum_{i=1}^{N}\sum_{l=1}^{N-1}\langle u_l(x_i), x_i-x_{i+l} \rangle \geq 0.\end{equation}
\end{definition}
Note that if each $u_{\ell }$ is $N$-cyclically monotone, then the family $(u_{1},u_{2},...,u_{N-1})$ is jointly $N$-monotone. Actually, one needs much less, since the $(N-1)$-tuplet $(u,u,...,u)$ is jointly $N$-monotone if and only if $u$ is $2$-monotone. On the other hand, $(u,0,0,...,0)$ is jointly $N$-monotone if and only if $u$ is $N$-monotone. See \cite{GG} for a complete discussion.

We now state the recent result of Galichon-Ghoussoub \cite{GG}, which establishes the remarkable duality between $N$-cyclically monotone operators, $N$-antisymmetric Hamiltonians and measure preserving $N$-involutions.
We shall also need the notion of  an \emph{$N$-sub-antisymmetric Hamiltonian on $\Omega$,} which is any function $H$ satisfying
\begin{equation}\label{subanti}
\hbox{$\sum\limits_{i=0}^{N-1}H(\sigma ^{i}(x_{1},...,x_{N}))\leq 0 $ on $%
\Omega^N$ and $H(x,x,..., x)=0$ for all $x\in \Omega$.}
\end{equation}
\begin{theorem} {\rm (Galichon-Ghoussoub)}
\label{var} Let $u_1,..., u_{N-1}:\Omega \to {\mathbb{R}}^d$ be bounded
measurable vector fields. The following properties are then equivalent:

\begin{enumerate}
\item The family $(u_1,..., u_{N-1})$ is jointly $N$-monotone a.e.,
that is there exists a measure zero set $\Omega_0$ such that $(u_1,...,
u_{N-1})$ is jointly $N$-monotone on $\Omega \setminus \Omega_0$.

\item The family $(u_1,..., u_{N-1})$ is in the polar of ${\mathcal{S}}_N(\Omega, \mu)$
in the following sense,
\begin{equation}\label{opti}
\inf\left\{\int_\Omega \sum\limits_{\ell=1}^{N-1}\langle u_\ell (x),
x-S^\ell x\rangle \, d\mu; S \in {\mathcal{S}}_N(\Omega, \mu)\right\}=0.
\end{equation}
\item There exists a $N$-sub-antisymmetric Hamiltonian $H$ which
is concave in the first variable, convex in the last $(N-1)$ variables such
that
\begin{equation}
\hbox{$(u_1(x),...,u_{N-1}(x))=\nabla_{2,..., N}H(x,x,...,x)$ \quad for a.e.
$x\in \Omega$.}
\end{equation}%
Moreover, $H$ is $N$-cyclically antisymmetric in the following sense: For a.e.  ${\bf x}=(x_1,..., x_{N})\in \Omega^N$, we have
\[
H(x_1, x_2, ...,x_N)+H_{2,..., N}(x_1, x_2,..., x_N)=0
\]
where $H_{2,..., N}$ is the concavification of the function $K({\bf x})=\sum\limits_{i=1}^{N-1}H(\sigma^i({\bf x}))$ with respect to the last $N-1$ variables.
\end{enumerate}
\end{theorem}
Note that (\ref{opti}) shows that the above is also equivalent to the statement that
\begin{equation}
\sup\{\int_{\Omega^N}\sum_{\ell=1}^{N-1}\langle u_\ell(x_1), x_{\ell +1}\rangle  d\pi({\bf x});\, \pi\in \P_{\rm sym}(\Omega^N, \mu)\}=\int_{\Omega}\sum_{\ell=1}^{N-1}\langle u_\ell(x), x \rangle\, d\mu (x),
\end{equation}
and that the supremum is attained at the image of $\mu$ by the map $x\to (x,x,...,x)$, which is nothing but a particular case of the symmetric Monge-Kantorovich problem, when the cost function is the one we are considering in (\ref{cost}) and when the family $(u_1,..., u_{N-1})$ is $N$-monotone.

Theorem \ref{main} below can now be seen as the extension of the above, when one considers an arbitrary family of $(N-1)$ vector fields. Indeed, note that in the case of the cost function (\ref{cost}),
\[
\ell_H(x)=\sup\left\{\langle u_1(x), x_2\rangle +....\langle u_{N-1}(x), x_N\rangle- H(x, x_2,..., x_N);\, (x_{2},...,x_{N}) \in \Omega^{N-1}\right\},
\]
which means that $\ell_H $ is essentially the standard Lagrangian associated to $H$ (i.e., Legendre transform of $H$ with respect to the last $N-1$-variables) and
 \[
\ell_H (x)=L_H(x, u_1(x), u_2(x),..., u_{N-1}(x)),
\]
where for $(x, p_1,...,p_{N-1})\in (\R^d)^N$,
\[
L_H(x, p_1,...,p_{N-1})=\sup\{\sum_{i=1}^{N-1}\langle p_i, y_i\rangle -H(x, y_1,..., y_{N-1}); y_i\in \Omega\}.
\]
The following result will be established in sections 3 and 4.

\begin{theorem}\label{main}
Given $(N-1)$ bounded vector fields $u_1, u_2,...., u_{N-1}$ from $\Omega$ to $\R^N$, and a probability measure $\mu$ on $\Omega$ that is absolutely continuous with respect to Lebesgue measure,
we consider the following variational problems:
\begin{eqnarray}
{\rm MK}_{\rm sym}:&=&\sup\{\int_{\Omega^N}\left[\langle u_1(x_1), x_2\rangle +....\langle u_{N-1}(x_1), x_N\rangle \right]  d\pi; \pi\in \P_{sym}(\Omega^N,\mu) \}.\label{dual}\\
{\rm DK}_{\rm sym}:&=&\inf \{\int_\Omega L_H(x, u_1(x), u_2(x),..., u_N(x))\, d\mu (x); H \in \H_N(\Omega)\}.\label{primal}\\
 {\rm MK}_{\rm cyc}:&=&\sup \{\int_{\Omega^N}\left[\langle u_1(x), Sx\rangle +\langle u_2(x), S^2x\rangle+....\langle u_{N-1}(x), S^{N-1}x\rangle \right] d\mu;\, S\in {\mathcal S}_N(\Omega) \}.
\end{eqnarray}
If ${\rm meas}(\partial \Omega)=0$, then the following holds:
\begin{enumerate}
\item ${\rm MK}_{\rm sym}={\rm DK}_{\rm sym}={\rm MK}_{\rm cyc}$.
\item  ${\rm MK}_{\rm cyc}$ is attained at some $S\in {\mathcal S}_N(\Omega)$, which means that ${\rm MK}_{\rm sym}$ is attained at an invariant measure $\pi_S$ that is the image of $\mu$ by the map $x\to (x, Sx, S^2x,..., S^{N-1}x)$.

\item There exists a function $H$ on $\R^{dN}$ that is
concave in the first variable, convex in the last $(N-1)$ variables and $N$-sub-antisymmetric on $\Omega$,
 such that
 \begin{equation}\label{inter1}
(u_1(x),..., u_{N-1}(x))\in \partial_{2,...,N} H (x, Sx, ..., S^{N-1}x)\quad a.e. \,\, x \in \Omega.
\end{equation}
Moreover, if either $u_i\in W^{1,1}_{loc} (\Omega)$ for $i=1,2,...,N-1$ or if $S$ is differentiable a.e., then there exists a $N$-cyclically antisymmetric Hamiltonian $H\in \H_N(\Omega)$ such that
\begin{equation}\label{rep}
(u_1(x),..., u_{N-1}(x))=\nabla_{2,...,N} H(x, Sx, ..., S^{N-1}x)\quad a.e. \,\, x \in \Omega.
\end{equation}
\item Assume that for any two families of points $x_1,..., x_N$ and $y_1,..., y_N$ in $\Omega$,
 the function
\[x\to\sum_{i=1}^{N-1}\langle u_i(x), y_i-x_i\rangle +\sum_{i=1}^{N-1}\langle u_i(y_{N-i})-u_i(x_{N-i}), x\rangle\]
has no critical point unless when $x_1=y_1.$ Then there exists a unique measure preserving $N$-involution $S$ such that $(\ref{rep})$ holds  for some concave-convex $N$-sub-antisymmetric Hamiltonian $H$.
 \end{enumerate}
\end{theorem}
If $u:\Omega \to \R^d$ is a single bounded vector field, then the above theorem applied to the family $(0,..., 0, u)$ yields the decomposition
\begin{equation}\label{rep.single}
(-u(Sx),0,...,0,u (x))=\nabla H(x, Sx, ..., S^{N-1}x)\qquad a.e. \quad x \in \Omega.
\end{equation}
If $S$ is the identity in the above representation,  $u$ is then $N$-cyclically monotone, which means that the above theorem essentially says that any bounded vector field is $N$-cyclically monotone up to a measure preserving $N$- involution. This is clearly in the same spirit as Brenier's theorem stating that any non-degenerate vector field is the gradient of a convex function (i.e., is $N$-cyclically monotone for all $N$) modulo a measure preserving transformation. Note that the representation of $2$-monotone operators as partial gradients of antisymmetric saddle functions was established by Krause \cite{Kra}. The general version of this result was established in \cite{GM} where it is shown that any bounded vector field is $2$-monotone up to a measure preserving involution. Theorem \ref{main} can be seen as an extension of this result  to the case where $N\geq 2$ and where there is more than one vector field.

Actually, in the case of a single vector field $u:\Omega \to \R^d$, one need not consider Hamiltonians on $\Omega^N$ as long as the requirement of $N$-antisymmetry is replaced by the following property:  Say that a function $F$ on $\R^d \times \R^d$ is \emph{$N$-cyclically sub-antisymmetric on $\Omega$}, if
\begin{equation}
\hbox{$F(x,x)=0$ and $\sum\limits_{i=1}^{N} F(x_i, x_{i+1})\leq 0$ for all cyclic families $x_1,..., x_N, x_{N+1}=x_1$ in $\Omega$.}
\end{equation}
Note that if a function $H(x_1,..., x_N)$ is $N$-sub-antisymmetric and if it only depends on the first two variables, then the function $F(x_1, x_2):=H(x_1, x_2,..., x_N)$
is $N$-cyclically sub-antisymmetric.

Our proof then yields the following result.

\begin{theorem}\label{main.2} Consider a vector field $u\in L^\infty(\Omega, \R^d)$,  then:
\begin{enumerate}
\item  For every $N\geq 2$, there exists a measure preserving $N$-involution $S$ on $\Omega$ and a globally Lipschitz concave-convex  function $F$ of $\R^d \times \R^d$ that is  $N$-cyclically sub-antisymmetric on $\Omega$, such that
\begin{equation}
\hbox{$(-u(Sx), u(x))\in \partial F(x,Sx)$ for a.e. $x\in \Omega$,}
\end{equation}
where $\partial H$ is the subdifferential of $H$ as a concave-convex
function \cite{Rock}.

\item If either $u\in W^{1,1}_{loc} (\Omega)$ or if $S$ is differentiable a.e., then
\begin{equation}\label{rep2}
\hbox{$u(x)=\nabla_2 F(x,Sx)$ for a.e. $x\in \Omega$.}
\end{equation}

\item Moreover $u$ is $N$-cyclically monotone on $\Omega$ if and only if $S=I$ in the representation (\ref{rep2}).

\end{enumerate}

\end{theorem}
Note that we cannot expect to have a function $F$ such that $\sum\limits_{i=1}^{N} F (x_i, x_{i+1})=0$ for all cyclic families $x_1,..., x_N, x_{N+1}=x_1$
in $\Omega$.
This is the reason why one needs to consider functions of $N$-variables in order to get $N$-antisymmetry as opposed to sub-antisymmetry. Note that the function defined by
\begin{equation}
H(x_1,x_2,..., x_N): =\frac{(N-1)F(x_1,x_2)- \sum_{i=2}^{N-1}F(x_i, x_{i+1})}{N},
\end{equation}
is $N$-antisymmetric in the sense of belonging to ${\mathcal H}_N(\Omega)$ while $H(x_1,x_2..., x_N) \geq F(x_1,x_2)$ on $\Omega^N$.

\section{The case of a general cost function}
Let $\mu_1, \mu_2, ..., \mu_N$ be a probability measure on a domain $\Omega$ of $\R^d$, and consider the following Monge-Kantorovich problem associated to a given cost function $c: \Omega^N\to \R\cup \{-\infty\}$.
\begin{equation}\label{gm1}
\hbox{${\rm MK}(c, \mu_1,..., \mu_N)=\sup\{\int_{\Omega^N}c(x_1,x_2,..., x_N) d\pi;\, \pi \in {\cal P}(\Omega^N)\, \& \, {\rm proj}_i\pi=\mu_i$ for all $i=1,...,N\}$, }
\end{equation}
where ${\cal P}(\Omega^N)$ is the set of probability measures on $\Omega^N$.
The following proposition is standard.

\begin{proposition} \label{uno}Assume $c$ is a finitely valued upper semi-continuous and bounded above cost function on $\Omega^N$, then
\begin{enumerate}
\item There exists $\pi_0\in \P (\Omega^N)$ with ${\rm proj}_i\pi_0=\mu_i$ for all $i=1,..., N$ where ${\rm MK}(c, \mu_1,..., \mu_N)$ is attained.

\item The following duality holds: ${\rm MK}(c, \mu_1,..., \mu_N)$ is equal to
\[
{\rm DK}(c, \mu_1,..., \mu_N):=\inf \biggl\{\int_{\Omega^N}\sum_{j=1}^N u_j(x_j)\, d\mu_j; \, (u_j)_{j=1}^N:\Omega \to \R^N\, {\rm \&}\,  \sum_{j=1}^Nu_j(x_j)\geq c(x_1,\dots,x_N)\biggr\},
\]
and there exists bounded borel functions $u^0_1, ..., u_N^0$ where ${\rm DK}(c, \mu_1,..., \mu_N)$ is attained.
\end{enumerate}
\end{proposition}
We now consider such Monge-Kantorovich problems in the presence of symmetry. Say that {\it $c$ is cyclically symmetric} if
\begin{equation}
\hbox{$c(x_1, x_2,..., x_N)=c(x_2, x_3,...,x_N, x_1)$ on $\Omega^N$.}
\end{equation}
\begin{proposition} \label{duo}Assume $c$ is a cyclically symmetric upper semi-continuous and bounded above cost function, and that all marginals $\mu_i$ are equal to $\mu$. Then
\begin{enumerate}
\item There exists $\tilde \pi_0\in \P_{\rm sym}(\Omega^N)$ with ${\rm proj}_i\pi_0=\mu$ for all $i=1,...,N$ where ${\rm MK}(c, \mu_1,..., \mu_N)$ is attained. Moreover,
\begin{equation}\label{gm2}
\hbox{${\rm MK(c,\mu_1,..., \mu_N)=MK}_{\rm sym}(c, \mu)=\sup\{\int_{\Omega^N}c(x_1,x_2,..., x_N) d\pi;\, \pi\in \P_{\rm sym}(\Omega^N)\, \, \&\, \, {\rm proj}_1\pi=\mu\}.$}
\end{equation}
\item If $c$ is finite, then
\begin{equation}\label{gm3}
\hbox{${\rm MK_{\rm sym}(c, \mu)=DK}_{\rm sym}^1(c, \mu):=\inf \biggl\{N\int_\Omega u(x)\, d\mu; \, u\in C_b(\Omega)$ with $ \sum_{j=1}^Nu(x_j)\geq c(x_1,\dots,x_N)\biggr\},$}
\end{equation}
and there exists a Borel and bounded function $u_0:\Omega \to \R$ so that
\begin{equation}
{\rm DK}_{\rm sym}^1(c, \mu)=N\int_\Omega u_0(x)\, d\mu.
\end{equation}
 Moreover, $u_0$ can be chosen in such a way that
\begin{equation}\label{iterate}
u_0(x)=\sup\left\{c(x, y_1, y_2, ..., y_{N-1})-\sum_{i=1}^{N-1}u_0(y_i)\right\}.
\end{equation}
\end{enumerate}

\end{proposition}
\noindent{\bf Proof:}  1. By Proposition \ref{uno}, ${\rm MK}(c, \mu,..., \mu)$ is attained at some $\pi_0\in {\cal P}(\Omega^N)$ with marginals ${\rm proj}_i\pi=\mu$ for $i=1,..., N$. Consider now the probability
$
\tilde \pi_0:=\frac{1}{N}\sum\limits_{i=1}^N\sigma^i\#\pi_0.
$
It is clearly $\sigma$-invariant and with marginal $\mu$. Since $c$ is cyclically symmetric, we have $\int_{\Omega^N}c(x_1,x_2,..., x_N) d\pi_0=\int_{\Omega^N}c(x_1,x_2,..., x_N) d\tilde \pi_0$, meaning that $\tilde \pi_0$ is also maximizing for both ${\rm MK}(c, \mu,..., \mu)$ and {\rm MK}$_{\rm sym}(c, \mu)$.\\

2. \, By Proposition \ref{uno}, ${\rm DK}(c, \mu)$ is attained at a family of bounded Borel functions $(u^0_i)_{i=1}^N$. Set
\[
u(x)=\frac{u^0_1(x)+u^0_2(x)+...+u^0_N(x)}{N},
\]
and note  that $c$ is cyclically symmetric,  $\sum_{j=1}^Nu(x_j)\geq c(x_1,\dots,x_N)$ and
\[
{\rm DK}(c, \mu)\leq {\rm DK}_{\rm sym}^1(c, \mu)\leq N\int_\Omega u(x)\, d\mu=\int_{\Omega^N}\sum_{j=1}^N u(x_j)\, d\mu=\int_{\Omega^N}\sum_{j=1}^N u^0_j(x_j)\, d\mu={\rm DK}(c; \mu,...,\mu).
\]
In order to show (\ref{iterate}) we consider the function
\[
\bar u(x)=\sup\left\{c(x, y_1, y_2, ..., y_{N-1})-\sum_{i=1}^{N-1}u(y_i);\, \, (y_1,..., y_{N-1}) \in \Omega^{N-1}\right\}.
\]
Since $u$ satisfies the constraint, we have $\bar u (x)\leq u (x)$.  We now claim that for all $(x, x_2,..., x_N)$, we have
\begin{equation}\label{iterate.1}
\bar u(x)+\sum_{i=2}^{N}\bar u(x_i)\geq Nc(x, x_2,..., x_{N})-(N-1)\big(\sum_{i=2}^{N}u(x_i)-u(x)\big).
\end{equation}
Indeed, by picking $y_i=x_{i+1}$ in the definition of $\bar u(x)$, we get that
\[
\bar u(x)\geq c(x, x_2, x_3,..., x_N)-\sum\limits_{i=2}^{N}u(x_i).
\]
Similary, pick $y_1=x$, $y_i=x_{i+1}$ in the definition of $\bar u(x_2)$,
\[
\bar u(x_2)\geq c(x_2, x, x_3,..., x_N)-u(x)-\sum\limits_{i=2}^{N}u(x_i) +u(x_2).
\]
Similary, pick $y_1=x$, $y_i=x_{i+1}$ in the definition of $\bar u(x_2)$,
\[
\bar u(x_3)\geq c(x_3, x_2, x,..., x_N)-u(x)-\sum\limits_{i=2}^{N}u(x_i)+u(x_3)
\]
Now add up the N above inequalities and use the fact that $c$ is symmetric to obtain (\ref{iterate}).
Consider now the function  $v (x)=\frac{\bar u(x)+(N-1)u(x)}{N}$ in such a way that $\bar u (x)\leq v (x) \leq u (x)$. Estimate  (\ref{iterate}) gives that
$
\sum_{i=1}^{N}v(x_i)\geq c(x_1, x_2, x_3,..., x_N),
$
hence $v$ satisfies the constraint in DK$_{\rm sym}^1(c, \mu)$. Let now
\[
u_0(x)=\inf\large\{w(x); \bar u \leq w \leq u\,  \& \, \sum_{i=1}^{N}w(x_i)\geq c(x_1, x_2, x_3,..., x_N)\large\}.
\]
The function $u_0$ clearly satisfies the constraint in DK$_{\rm sym}^1(c,\mu)$. Note also that
\begin{eqnarray*}
 \bar u_0(x)&=&\sup\big\{c(x, y_1, y_2, ..., y_{N-1})-\sum_{i=1}^{N-1}u_0(y_i);\, (y_1,..., y_{N-1}) \in \Omega^{N-1}\big\}\\
&\geq &\sup\big\{c(x, y_1, y_2, ..., y_{N-1})-\sum_{i=1}^{N-1}u(y_i);\, (y_1,..., y_{N-1}) \in \Omega^{N-1}\big\}\\
&=&\bar u(x),
\end{eqnarray*}
which means that $u\geq \bar u_0\geq \bar u$, If now $\bar u_0(\bar x)<u_0(\bar x)$ for some $\bar x\in \Omega$,  then $u_0(\bar x)>\bar u_0(\bar x)\geq \bar u(\bar x)$, hence contradicting the minimality of $u_0$. It follows that $u_0=\bar u_0$, and since $u_0\leq u$, it does minimize the functional in DK$_{\rm sym}^1(c, \mu)$.\hfill $\Box$\\

We now consider the symmetric Monge-Kantorovich problem when $c$ is not assumed to be cyclically symmetric.\\

\noindent {\bf Proof of Theorem \ref{symmetric.dual}:} let $\tilde c$ is the symmetrized of $c$ defined for any ${\bf x}=(x_1, x_2,..., x_N)$ by
$
\tilde c({\bf x})=\frac{1}{N}\sum_{i=1}^N c(\sigma^i{\bf x}).
$
It is clear that ${\rm MK}_{\rm sym}(c, \mu)={\rm MK}_{\rm sym}(\tilde c, \mu)$, which by the last proposition is equal to ${\rm DK}_{\rm sym}^1(\tilde c, \mu)$.

On the other hand, For $H\in {\mathcal H}_N(\Omega)$, we let $\ell_H^c$ be the $c$-Legendre transform of $H$ with respect to the last $(N-1)$ variables, that is
\[
\ell_H^c(x)=\sup\left\{c(x, x_2,..., x_N)-H\left( x,x_2,..., x_N\right); (x_{2},...,x_{N})\in \Omega^{N-1}\right\}
\]
For any invariant probability measure $\pi$, with $1$-marginal $\mu$, any $\ell \in L^1(\Omega; \mu)$ and any $H\in \H$, we have
\begin{equation*}
\int_{\Omega^N}c(x_1,x_2,..., x_N) d\pi = \int_{\Omega^N} \left[c(x_1,x_2,..., x_N) -\ell\left( x_{1}\right) - H(x_1,x_2,..., x_N)\right]d\pi \left( {\bf x}\right) +\int_\Omega \ell (x_1) d\mu (x_1),
\end{equation*}%
hence if $\ell(x_1) \geq c(x_1, ..., x_N)-H (x_1,..., x_N)$ for all $(x_1, x_2,..., x_N)$, then
$
\int_{\Omega^N}c(x_1,x_2,..., x_N) d\pi \leq \int_\Omega \ell (x) d\mu (x)$,
and therefore ${\rm MK}_{\rm sym}(c, \mu)\leq  {\rm DK}^2_{\rm sym}(c, \mu)$.

For the reverse inequality, we shall use the fact that
\[
{\rm MK}_{\rm sym}(c, \mu)={\rm MK}_{\rm sym}(\tilde c, \mu)={\rm DK}^1_{\rm sym}(\tilde c, \mu)=N\int_\Omega u_0(x) d\mu,
\]
where $u_0$ is a lower semi-continuous function satisfying
\[
u_0(x)=\sup\Big\{\tilde c(x,y_1\dots,y_{N-1})-\sum_{i=1}^{N-1}u_0(y_i):\, y_1,..., y_{N-1}\in\R^d\Big\}.
\]
Notice that  $Nu_0=\ell_H^{\tilde c}$, where $H$ is the $N$-cyclically symmetric Hamiltonian  defined by
\[H(x_1,x_2,...,x_N):= \sum_{i=2}^{N}u_0(x_i)-(N-1)u_0(x_1).\]
Finally, it is easy to check that
$Nu_0=\ell_{H_\infty}^{c}$, where
$H_\infty$ is the $N$-cyclically symmetric Hamiltonian
\[
H_\infty(x_1,x_2,...,x_N):= -\frac{1}{N}\sum_{i=1}^N c(\sigma^i(x_1, x_2,..., x_N)) +c(x_1, x_2,..., x_N)+\sum_{i=2}^{N}u_0(x_i)-(N-1)u_0(x_1).
\]
It follows that
\[
{\rm MK}_{\rm sym}(c, \mu)=N\int_\Omega u_0(x) d\mu={\rm DK}^2_{\rm sym}(c, \mu),
\]
and that the latter is attained at $H_\infty$. This completes the proof of Theorem \ref{symmetric.dual}
 \hfill $\square$\\
\noindent {\bf Proof of Theorem \ref{extremal}:} This will follow from Theorem \ref{symmetric.dual} combined with the following three lemmas.

\begin{lemma}  Let $L:=\ell^c_{H_\infty}$ where $H_\infty$ is a fixed Hamiltonian in $\H_N(\Omega)$, and let $x\in \Omega$ be such that $(x_2, x_3,..., x_N) \to c(x, x_2,..., x_N)-H_\infty(x,x_2,..., x_N)$ attains its maximum uniquely at $S_1x,S_2x, ..., S_{N-1}x$.
 Let $H\in \H$, $r\in \R$ and consider  $L^c_r:=\ell^c_{H_\infty +rH}$ to be the $c$-Legendre transform associated to the Hamiltonian $H_\infty +rH$. Then, we have
\begin{eqnarray}
\lim_{r \to 0}\frac{L_r(x)-L(x)}{r}=H(x, S_1x, S_2x,...,S_{N-1}x).
\end{eqnarray}
\end{lemma}
{\bf Proof:} Let $(x_{1,r},...., x_{N-1, r})$ be points in $\Omega$, where
\[
(x_2, x_3,..., x_N) \to c(x, x_2,..., x_N)-H_\infty(x,x_2,..., x_N)-rH (x,x_2,..., x_N)
\]
attains its maximum. It follows that
\begin{eqnarray*}
L_r(x)-L(x)&=&c(x, x_{1,r},...., x_{N-1, r})-H_\infty (x, x_{1,r},...., x_{N-1, r})-rH(x, x_{1,r},...., x_{N-1, r})\\
&&-c(x, S_1x, S_2x,..., S_{N-1}x)+H_\infty(x, S_1x,S_2x, ..., S_{N-1}x),
\end{eqnarray*}
and therefore
\begin{eqnarray*}
-H_\infty(x, S_1x,S_2x, ..., S_{N-1}x)\leq \frac{L_r(x)-L(x)}{r}\leq -H(x, x_{1,r},...., x_{N-1, r}).
\end{eqnarray*}
Since $S_1x, ..., S_2x$ are unique maxima, it follows that as $r\to 0$, we have that $x_{i,r}$ converges to $S_ix$, from which we conclude that $\lim_{r \to 0}\frac{L_r(x)-L(x)}{r}=H(x, S_1x, S_2x,...,S_{N-1}x)$.

\begin{lemma} Assume that
\[
{\rm DK}^2_{\rm sym}(c, \mu):=\inf\Large\{\int_\Omega \ell_H^c(x)d\mu (x); H\in \H_N(\Omega)\Large\}
\]
 is attained at some $H_\infty \in \H_N(\Omega)$, and that for $\mu$-almost $x\in \Omega$ the map
\[
(x_2, x_3,..., x_N) \to c(x, x_2,..., x_N)-H_\infty(x,x_2,..., x_N)
\]
 attains its maximum uniquely at $S_1x,S_2x, ..., S_{N-1}x$.   Then, for any $H \in \H_N(\Omega)$, we have
\[
\int_\Omega H (x, S_1x,S_2x, ..., S_{N-1}x)\, d\mu=0.
\]
\end{lemma}

\noindent{\bf Proof:} Let $L=\ell^c_{H_\infty}$. For any $N$-symmetric Hamiltonian $H$ and $r\in \R$,  consider  $L^c_r=\ell^c_{H_\infty +rH}$ to be the $c$-Legendre transform associated to the Hamiltonian $H_\infty +rH$. The above lemma yields that for $\mu$-almost all $x\in \Omega$, we have
\begin{eqnarray*}
\lim_{r \to 0}\frac{L_r(x)-L(x)}{r}=H(x, S_1x, S_2x,...,S_{N-1}x).
\end{eqnarray*}
On the other hand, the extremality of $L:=\ell_{H_\infty}$ gives that
\begin{eqnarray*}
0=\lim_{r \to 0}\int_\Omega \frac{L_r(x)-L(x)}{r}d\mu=\int_\Omega H(x, S_1x, S_2x,...,S_{N-1}x).
\end{eqnarray*}\hfill $\Box$

At the core  of our results is the following duality between $N$-antisymmetric Hamiltonians and measure preserving $N$-involutions, which will be crucial to what follows.

\begin{lemma} \label{magic} Let $S_1, S_2, ..., S_{N-1}$ be $\mu$-measurable maps on $\Omega$. The following statements are then equivalent:

\begin{enumerate}
\item $\int_\Omega H (x, S_1x,S_2x, ..., S_{N-1}x) d\mu =0$ for any $N$-cyclically symmetric Hamiltonian $H$.

\item There exists $S:\Omega \to \Omega$, $\mu$-measure preserving such that
  $S^N=I$ and $S^i=S_i$ for all $i=1, ..., N-1$.
\end{enumerate}
\end{lemma}
\noindent{\bf Proof:}  If $S$ is $\mu$-measure preserving and $S^N=I$ a.e.,  then
\[
\int_\Omega H(x, Sx, S^2x,...S^{N-1}x)d\mu=\int_\Omega H(Sx, S^2x,...S^{N-1}x, x)d\mu=...=\int_\Omega H(S^{N-1}x,..., S^2x, x, Sx)d\mu
\]
Since $H$ is $N$-symmetric, then
\[
H(x, Sx, S^2x,...S^{N-1}x)+H(Sx, S^2x,...S^{N-1}x, x)+...H(S^{N-1}x,..., S^2x, x, Sx)=0.
\]
It follows that $\int_\Omega H(x, Sx, S^2x,...S^{N-1}x)d\mu=0$.\\

For the reverse implication, assume $\int_\Omega H (x, S_1x,S_2x, ..., S_{N-1}x) d\mu=0$ for any $N$-cyclically symmetric Hamiltonian $H$.   By using the identity with Hamiltonians $(H_i)_{i=1}^N$ of the form
\[
H_i(x_1, x_2,..., x_N)=f(x_1)-f(x_i),
\]
where $f$ is a continuous function, one gets that $S_i$ is measure preserving for each $i=1,..., N-1$.

Now take for each fixed  $i=1, ..., N$, the Hamiltonian
\[
H_i(x_1, x_2,..., x_N)=|x_i-S_1^ix_N|-|S_1^ix_1-x_{i+1}|-|x_{i+1}-S_1^ix_1|+|S_1^ix_2-x_{i+2}|.
\]
We have that $H_i\in \H_N(\Omega)$ for each $i$, since it is of the form $H_i(x_1,..., x_N)=f(x_1, x_i, x_N)-f(x_2, x_{i+1}, x_1)$. Hence, for each $i$,
\[
0=\int_\Omega H_i (x, S_1x,S_2x, ..., S_{N-1}x) d\mu=0=\int_\Omega ((|S_{i-1}x-S_1^iS_{N-1}|+|S_1^iS_1x-S_{i+1}x|)d\mu=0.
\]
It follows that $S_{i+1}=S_1^{i+1}$ and $S_{i-1}x=S_1^iS_{N-1}$ for each $i=1,...,N$. The latter applied to $i=1$, yields $x=S_1S_{N-1}=S_1S_1^{N-1}=S_1^N$, and we are done.

\section{Concave-convexification of $N$-antisymmetric functions}

 Let $\Omega$ be a bounded domain in $\R^d$, and consider  the class
\begin{eqnarray}
{\cal H}^{^-}_N(\Omega):=\big \{H \in C(\bar \Omega^N);\,  \sum\limits_{i=0}^{N-1}H (\sigma^i({\bf x}))\leq 0  \, \text{ for all} \, {\bf x} \in \Omega^N\}.
\end{eqnarray}
For each $H \in {\cal H}^{^-}_N(\Omega)$, we associate the following functional on $\Omega \times (\R^d)^{N-1}$,
\begin{equation}\label{legendre*}
L_H(x, p_1,...,p_{N-1})=\sup\left\{\sum_{i=1}^{N-1}\langle p_i, y_i\rangle -H(x, y_1,..., y_{N-1}); y_i\in \Omega\right\}.
\end{equation}
Denote by
\[
{\cal L}_-(N)=\{L_H;H \in {\cal H}^{^-}_N(\Omega) \}.
\]
  Our plan is to show that one can associate to $H,$
 \begin{itemize}
 \item a globally Lipschitz-continuous function $H^1_{reg} \in {\cal L}_-(N)$ that is concave in the first variable, convex in the last $(N-1)$
variables such that  $L_{H^1_{reg}} \leq L_{H}.$

 \item a globally Lipschitz-continuous  function $H^2_{reg} \in {\cal L}(N)$ such that $H_{reg}^2\geq H_{reg}^1$ and hence
\[
L_{H^2_{reg}}\leq L_{H^1_{reg}} \leq L_{H}.
\]

 \end{itemize}
 Suppose that $\Omega$ is contained in a ball $B_{R}$ centered at the origin with radius $R>0$ in $\R^d$, we define {\it ``an $(\bar \Omega \times B_R)$ restricted Legendre transform"} of $L_H$ as
\[L^*_H(p_1,..., p_{N-1},,x) =\sup_{q\in \bar \Omega, y_i \in B_R}\left\{
\langle q,x \rangle+ \sum_{i=1}^{N-1}\langle p_i, y_i\rangle - L_{H}(q, y_1, y_2,..., y_{N-1})\right\}.
\]
Similarly, we define on $\R^d \times (\R^d)^{N-1}$,
\begin{eqnarray}
L_{H}^{**}(x,p_1,..., p_{N-1})=\sup_{p\in \bar \Omega, x_i\in B_R}\left\{\langle x,p\rangle+
 \sum_{i=1}^{N-1}\langle p_i, x_i\rangle- L^*_{H}(x_1,..., x_{N-1},p)\right\}.
\end{eqnarray}
For any function  $L:\R^d \times (\R^d)^{N-1} \to \R$, we define its {\it ``$B_{R}$-Hamiltonian"} by
\begin{eqnarray}
H_{L}(x,y_1,..., y_{N-1})= \sup_{ p_i \in B_R}\left\{\sum_{i=1}^{N-1}\langle p_i, y_i\rangle- L(x,p_1,..., p_{N-1})\right\}.
\end{eqnarray}
Finally, for $H \in {\cal H}^{^-}_N(\Omega)$, we consider the following two regularizations of $H$:
\begin{equation}
H^1_{reg}({\bf x})=H_{L_H^{**}}({\bf x}),
\end{equation}
and
\begin{equation}
H^2_{reg}({\bf x})=\frac{(N-1)H^1_{reg}({\bf x})-\sum_{i=1}^{N-1}H^1_{reg}(\sigma^{i}({\bf x}))}{N}.
\end{equation}
We list some of the properties of $H^1_{reg}$, $H^2_{reg}$, $L^1_{ H_{reg}}$ and $L^1_{ H_{reg}}$.

\begin{proposition} \label{prop1} If $H \in{\cal H}^{^-}_N(\Omega)$, then the following statements hold:
\begin{enumerate}
\item   $ H^1_{reg}$ is a concave-convex on $\R^d\times \R^{d(N-1)}$ whose restriction to $\bar \Omega^N$ belong to ${\cal H}^{^-}_N(\Omega)$.
\item $ H^2_{reg}$ belongs to $ {\cal
H}_N(\Omega)$,  and $H^2_{reg} \geq  H^1_{reg}$ on $\bar {\Omega}^{N}$.
\item   $L_{ H^1_{reg}}$ is convex and continuous in all variables and  $L_{ H^2_{reg}}\leq L_{ H^1_{reg}} \leq L_{ H}$ on $\bar \Omega \times (B_R)^{N-1}$.
\item   $|L_{ H^1_{reg}}(x,p_1,..., p_{N-1})| \leq R\|x\|+R\sum\limits_{i=1}^{N-1}\|p_i\| + (2N+1)R^2$ for all $x$ and all $(p_i)_{i=1}^{N-1}$ in $\R^d$.
\item $| H^1_{reg}(x,y_1,..., y_{N-1})| \leq R\|x\|+R\sum\limits_{i=1}^{N-1}\|y_i\| + 2NR^2$ for all $x$ and all $(y_i)_{i=1}^{N-1}$ in $\R^d$.

\item    $L_{ H^2_{reg}}$ and $H^2_{reg}$ are both Lipschitz continuous with Lipschitz constants less than
$4NR.$
\end{enumerate}
\end{proposition}
The proof will require several lemmas.

\begin{lemma}\label{ine2} With the above notation, we have the following properties:   \begin{enumerate}

\item $L^{**}_H(x,p_1,..., p_{N-1}) \leq L_H(x,p_1,...,p_{N-1})$ for $x\in \bar \Omega$ and $p_i \in \R^d$ for $i=1,..., N-1$.

\item If $H^1_{\rm reg}$ denotes $H_{L^{**}_H}$, then  $H^1_{\rm reg}$ is concave  in the first variable and convex  in the last $(N-1)$ variables.
\item $L_{H^1_{\rm reg}}$ is jointly convex in all variables.
\end{enumerate}
\end{lemma}
\textbf{Proof.} 1) For $x\in \bar \Omega$ and $p_i \in \R^d$, $i=1,..., N-1$, we have
 we have
  \begin{eqnarray*}\label{l8}
  L_{H}^{**}(x,p_1,..., p_{N-1})&=&\sup_{q\in \bar \Omega, r_i\in B_R}\left\{\langle x,q\rangle+
 \sum_{i=1}^{N-1}\langle p_i, r_i\rangle- L^*_H(r_1,..., r_{N-1},q)\right\}\\
 &=&\sup_{q\in \bar \Omega, r_i\in B_R}\left\{\langle x,q\rangle+
 \sum_{i=1}^{N-1}\langle p_i, r_i\rangle-\sup_{y\in \bar \Omega, y_i \in B_R}\{
\langle y,q \rangle+ \sum_{i=1}^{N-1}\langle r_i, y_i\rangle - L_{H}(y, y_1,..., y_{N-1})\}\right\}\\
 &=&\sup_{q\in \bar \Omega, r_i\in B_R}\inf_{y\in \bar \Omega, y_i \in B_R}\left\{\langle x,q\rangle+
 \sum_{i=1}^{N-1}\langle p_i, r_i\rangle-\langle y, q\rangle- \sum_{i=1}^{N-1}\langle r_i, y_i\rangle +L_{H}(y, y_1,..., y_{N-1})\right\}\\
 &=&\sup_{q\in \bar \Omega, r_i\in B_R}\inf_{y\in \bar \Omega, y_i \in B_R}\left\{\langle q,x-y\rangle+
 \sum_{i=1}^{N-1}\langle p_i-y_i, r_i\rangle+L_{H}(y, y_1,..., y_{N-1})\right\}\\
 &=&\sup_{q\in \bar \Omega, r_i\in B_R}\inf_{y\in \bar \Omega, y_i \in B_R}\left\{\langle q,x-y\rangle+
 \sum_{i=1}^{N-1}\langle p_i-y_i, r_i\rangle+\sup\limits_{t_i\in \Omega}\{\sum_{i=1}^{N-1}\langle t_i, y_i\rangle -H(y, t_1,..., t_{N-1}) \}
\right\}\\
&=&\sup_{q\in \bar \Omega, r_i\in B_R}\inf_{y\in \bar \Omega, y_i \in B_R}\sup\limits_{t_i\in \Omega}\left\{\langle q,x-y\rangle+
 \sum_{i=1}^{N-1}\langle p_i-y_i, r_i\rangle+\sum_{i=1}^{N-1}\langle t_i, y_i\rangle -H(y, t_1,..., t_{N-1})
\right\}\\
&=&\inf_{y\in \bar \Omega, y_i \in B_R}\sup_{q\in \bar \Omega, r_i\in B_R}\sup\limits_{t_i\in \Omega}\left\{\langle q,x-y\rangle+
 \sum_{i=1}^{N-1}\langle p_i-y_i, r_i\rangle+\sum_{i=1}^{N-1}\langle t_i, y_i\rangle -H(y, t_1,..., t_{N-1})
\right\}.
 \end{eqnarray*}
 By taking $y=x$ and $y_i=p_i$, we readily get that $L^{**}_H(x,p_1,..., p_{N-1}) \leq L_H(x,p_1,...,p_{N-1})$.\\

 For 2) note first that by definition
 \[
H_{L^{**}}(x,y_1,..., y_{N-1})= \sup_{ p_i \in B_R}\left\{\sum_{i=1}^{N-1}\langle p_i, y_i\rangle- L_H^{**}(x,p_1,..., p_{N-1})\right\},
\]
and therefore for all $x \in \R^d,$ the function $(y_1,...,y_{N-1}) \to H_{L^{**}}(x,y_1,..., y_{N-1})$ is convex. We shall show that  for all $(y_1,...,y_{N-1})\in (\R^d)^{N-1}$,  the function
$x\to H_{L^{**}}(x,y_1,..., y_{N-1})$ is concave. In fact we show that
\[
x \to -H_{L^{**}}(x,y_1,..., y_{N-1})=\inf_{ p_i \in B_R}\{L_H^{**}(x,p_1,..., p_{N-1})-\sum_{i=1}^{N-1}\langle p_i, y_i\rangle\}
\]
is convex. Indeed, consider $ \lambda  \in (0,1)$ and elements $x_1, x_2 \in \R^d$, then for any $a, b$ such that
\[
\hbox{$a > -H_{L^{**}}(x_1,y_1,..., y_{N-1})$ and $b > -H_{L^{**}}(x_2,y_1,..., y_{N-1})$,}
\]
 we can find $(r_i)_{i=1}^{N-1}$ and  $(q_i)_{i=1}^{N-1}$  in $(\R^d)^{N-1}$ such that
\[-H_{L^{**}}(x_1,y_1,..., y_{N-1}) \leq L_H^{**}(x_1,r_1,..., r_{N-1})-\sum_{i=1}^{N-1}\langle r_i, y_i\rangle \leq a,
\]
and
\[ -H_{L^{**}}(x_2,y_1,...,y_{N-1}) \leq L_H^{**}(x_2,q_1,...,q_{N-1})-\sum_{i=1}^{N-1}\langle q_i, y_i\rangle \leq b.
\]
Use the convexity of the ball $B_R$ and the convexity of the function $L_H^{**}$ in both variables  to write
\begin{eqnarray*}
-H_{L_H^{**}}(\lambda x_1+(1-\lambda) x_2, y_1,...,y_{N-1}) &=& \inf_{ p_i \in B_R}\{L_H^{**}(\lambda x_1+(1-\lambda) x_2, p_1,...,p_{N-1} )-\sum_{i=1}^{N-1}\langle p_i, y_i\rangle\}\\ & \leq & L_H^{**}(\lambda x_1+(1-\lambda) x_2, \lambda r_1 +(1-\lambda)q_1,..., \lambda r_{N-1} +(1-\lambda)q_{N-1}))\\
&&-\sum_{i=1}^{N-1}\langle \lambda r_i +(1-\lambda)q_i, y_i\rangle\\
& \leq & \lambda \big (L_H^{**}( x_1, r_1,..., r_{N-1})-\sum_{i=1}^{N-1}\langle r_i,y_i\rangle \big ) \\
&&+(1-\lambda) \big (L_H^{**}( x_2, q_1,..., q_{N-1})-\sum_{i=1}^{N-1}\langle q_i,y_i\rangle \big ) \}\\
& \leq & \lambda a +(1-\lambda)b,
\end{eqnarray*}
which establishes the concavity of $x \to H_{L_H^{**}}(x,y_1,..., y_{N-1})$. It then follows that $L_{H^1_{reg}}=L_{H_{L_H^{**}}}$ is convex in all variables that proves part 3).

\begin{lemma}\label{ine1} If $H \in {\cal H}^{^-}_N(\Omega)$, then $H^1_{\rm reg} \in {\cal H}^{^-}_N(\Omega).$
\end{lemma}

\textbf{Proof.}  Let $i,j=1, 2,.., N$.   We first  show that
\begin{eqnarray}\label{wha}
 \sum_{i=1}^{N}\Big \{ \sum_{j=1, j\not =i}^{N}\langle p_j^i, x_j\rangle- L^{**}_H(R^{i-1}(p^i_1,..., p_{i-1}^i,x_i, p_{i+1}^i,..., p^i_{N}) ) \Big \}\leq 0,
\end{eqnarray}
for all  $x_i \in \Omega $ and $p_j^i \in \R^d.$ Indeed, we have
\begin{eqnarray*}\label{legendre*}
L_H(\sigma^{i-1}(p^i_1,..., p_{i-1}^i,x_i, p_{i+1}^i,..., p^i_{N}) )&=&\sup\left\{\sum_{j=1, j\not =i}^{N}\langle p_j^i, y_j\rangle -H(\sigma^{i-1}( y_1,..., y_{i-1}, x_i,  y_{i+1},..., y_{N})); y_j\in \Omega\right\}\\
&\geq& \sum_{j=1, j\not =i}^{N}\langle p_j^i, x_j\rangle -H(\sigma^{i-1}(x_1,x_2,..., x_n)).
\end{eqnarray*}
Taking summation over $i$ implies that

\begin{eqnarray*}
\sum_{i=1}^{N} L_H(\sigma^{i-1}(p^i_1,..., p_{i-1}^i,x_i, p_{i+1}^i,..., p^i_{N}) ) \geq \sum_{i=1}^{N}\sum_{j=1, j\not =i}^{N}\langle p_j^i, x_j\rangle-\sum_{i=1}^{N}H(\sigma^{i-1}(x_1,x_2,..., x_n))
\end{eqnarray*}
Since $\sum_{i=1}^{N}H(\sigma^{i-1}(x_1,x_2,..., x_n)) \leq 0,$ we obtain
\begin{eqnarray*}
\sum_{i=1}^{N} L_H(\sigma^{i-1}(p^i_1,..., p_{i-1}^i,x_i, p_{i+1}^i,..., p^i_{N}) )
&\geq& \sum_{i=1}^{N}\sum_{j=1, j\not =i}^{N}\langle p_j^i, x_j\rangle.
\end{eqnarray*}
It follows from the  definition of $L_H^{**}$ that

\begin{eqnarray*}
\sum_{i=1}^{N} L^{**}_H(\sigma^{i-1}(p^i_1,..., p_{i-1}^i,x_i, p_{i+1}^i,..., p^i_{N}) )\geq \sum_{i=1}^{N}\sum_{j=1, j\not =i}^{N}\langle p_j^i, x_j\rangle.
\end{eqnarray*}
By moving the left hand side expression to the the other side, we have
\begin{eqnarray*}
0\geq \sum_{i=1}^{N}\Big \{ \sum_{j=1, j\not =i}^{N}\langle p_j^i, x_j\rangle- L^{**}_H(\sigma^{i-1}(p^i_1,..., p_{i-1}^i,x_i, p_{i+1}^i,..., p^i_{N}) ) \Big \}.
\end{eqnarray*}
Taking sup over all $p_i^j \in B_R$ we obtain
$ \sum_{i=1}^{N}H_{L^{**}_H}(\sigma^{i-1}(x_1,x_2,..., x_n))\leq 0$ and we are done.
  \hfill $\square$\\

We now recall the following standard elementary result.
\begin{lemma}\label{Gp} Let $D$ be an open set in $\R^m$ such that $\bar D \subset \tilde B_R$ where $\tilde B_R$ is ball with radious $R$ centered at the origin in $\R^m.$ Let $f:\R^m \to \R$ and define $\tilde f : \R^m \to \R$  by
\[\tilde f(y)= \sup_{z \in D}\{\langle y,z\rangle-f(z)\}.\]
If $f \in L^{\infty}(D),$ then $\tilde f$ is a convex Lipschitz function and
\[
\hbox{$|\tilde f (y_1)-\tilde f(y_2)| \leq R \|y_1-y_2\|$ for all $y_1, y_2 \in \R^m.$}\]

\end{lemma}
\begin{lemma} \label{prop2} If $H \in {\cal H}^{^-}_N(\Omega)$, then the following statements hold:
\begin{enumerate}
\item  $|L^{**}_H(x,p_1,..., p_{N-1})| \leq R\|x\|+R\sum_{i=1}^{N-1}\|p_i\| + (2N-1)R^2$ for all $x$ and  $(p_i)_{i=1}^{N-1}$ in $\R^d$.
\item $|H_{L^{**}_H}(x,y_1,..., y_{N-1})| \leq R\|x\|+R\sum_{i=1}^{N-1}\|y_i\| + 2NR^2$ for all $x$ and $(y_i)_{i=1}^{N-1}$ in $\R^d$.
\item $L^{**}_{H}$ and  $H_{L^{**}_H}$ are Lipschitz continuous with Lipschitz constants $Lip(H_{L^{**}_H}), Lip(L^{**}_H)\leq NR.$
\end{enumerate}
\end{lemma}

\textbf{Proof.} Since $H$ is $N$-sub-antisymmetric, we have $H(x,..., x)\leq 0$,
hence
\[
\hbox{$L_H(x, p_1,...,p_{N-1})\geq \sum\limits_{i=1}^{N-1}\langle p_i, x\rangle$ \,\,  on $\bar \Omega \times (\R^d)^{N-1}$.}
\]
This together with the fact that $\bar \Omega \subset B_R$ imply that
\begin{eqnarray*}
L_{H}^*(p_1,...,p_{N-1},x)&=&\sup_{q\in \bar \Omega, y_i \in B_R}\left\{
\langle q,x \rangle+ \sum_{i=1}^{N-1}\langle p_i, y_i\rangle - L_{H}(q, y_1, y_2,..., y_{N-1})\right\}.
\\
&\leq &\sup_{q\in \bar \Omega, y_i \in B_R}\left\{
\langle q,x \rangle+ \sum_{i=1}^{N-1}\langle p_i, y_i\rangle - \sum\limits_{i=1}^{N-1}\langle q, y_i\rangle\right\}.
\\
&\leq&R\|x\|+R\sum_{i=1}^{N-1}\|p_i\| + (N-1)R^2.
\end{eqnarray*}
With a similar argument we obtain that $L^{**}_H(x,p_1,..., p_{N-1}) \leq R\|x\|+R\sum_{i=1}^{N-1}\|p_i\| + (N-1)R^2.$ We also have
\begin{eqnarray*}
L_{H}^{**}(x,p_1,..., p_{N-1})&=&\sup_{p\in \bar \Omega, x_i\in B_R}\left\{\langle x,p\rangle+
 \sum_{i=1}^{N-1}\langle p_i, x_i\rangle- L^*_{H}(x_1,..., x_{N-1},p)\right\}\\
&\geq & \langle x,p\rangle+
 \sum_{i=1}^{N-1}\langle p_i, x_i\rangle- L^*_{H}(x_1,..., x_{N-1},p)\\
 &\geq &- R\|x\|-R\sum_{i=1}^{N-1}\|p_i\| -R\|p\|-R\sum_{i=1}^{N-1}\|x_i\| -(N-1)R^2\\
 &\geq& - R\|x\|-R\sum_{i=1}^{N-1}\|p_i\| -(2N-1)R^2.
\end{eqnarray*}
Therefore $ |L^{**}_H(x,p_1,..., p_{N-1})| \leq R\|x\|+R\sum_{i=1}^{N-1}\|p_i\| + (2N-1)R^2.$  The estimate for $H_{L^{**}_H}$  can be easily deduced from its definition together with the  estimate on $L^{**}_H.$ This completes the proof of part (1).\\

For (2) set $D=  \Omega \times \Pi_{i=1}^{N-1}B_R$, then $D \subset \tilde B_{NR}$ where  $\tilde B_{NR}$ is a ball with radius $NR$ in  $\R^{dN}.$ Now assuming $f=L^*_H$ in Lemma \ref{Gp}, we have that $\tilde f =L^{**}_H$. Therefore $L^{**}_H$ is Lipschitz in $(\R^d)^N$ with $Lip(L^{**}_H) \leq NR.$ To prove that $H_{L^{**}_H}$ is Lipschitz continuous, we first fix $y \in \R^{d}$ and define $f_y: (\R^d)^{N-1} \to \R$ by
\[
f_y(p_1,...,p_{N-1})=L^{**}_H(y,p_1,...p_{N-1}).
\]
 Assuming $D=B_R \subset \R^{N}$ in Proposition \ref{Gp}, we obtain that the map
 \[
 (x_1,..., x_{N-1}) \to \tilde f_y(x_1,..., x_{N-1})=H_{L^{**}_H}(y, x_1,..., x_{N-1})
 \]
  is Lipschitz and
\begin{equation}\label{lip}
|H_{L^{**}_H}(y, x_1,..., x_{N-1})-H_{L^{**}_H}(y, z_1,..., z_{N-1})| \leq R \sum_{i=1}^{N-1}\|x_i-z_i\|
\end{equation}
for all $(x_i), (z_i) \in (\R^d)^{N-1}.$ Noticing that the Lipschitz constant $R$ is independent of $y,$ the above inequality holds for all $(x_i), (z_i)\in (\R^d)^{N-1}$ and $ y \in \R^{d}.$ To prove $H_{L^{**}_H}(y, x_1,..., x_{N-1})$ is Lipschitz with respect to the first variable $y$,  let  $r>0$ and $y_1, y_2 \in \R^{d}.$ Let $p_1,..., p_{N-1}$ and  $q_1,..., q_{N-1}$ be such that
\[
\sum_{i=1}^{N-1}\langle x_i,q_i\rangle- L^{**}_{H}(y_1,q_1,..., q_{N-1}) \leq H_{L^{**}_H}(y_1, x_1,..., x_{N-1}) \leq \sum_{i=1}^{N-1}\langle x_i,p_i\rangle- L^{**}_{H}(y_1,p_1,..., p_{N-1})+ r,\]
and
\[\sum_{i=1}^{N-1}\langle x_i,p_i\rangle- L^{**}_{H}(y_2,p_1,..., p_{N-1}) \leq H_{L^{**}_H}(y_2, x_1,..., x_{N-1}) \leq \sum_{i=1}^{N-1}\langle x_i,q_i\rangle- L^{**}_{H}(y_2,q_1,..., q_{N-1})+ r.\]
It follows that
\begin{eqnarray*}
L^{**}_{H}(y_2,q_1,..., q_{N-1})-L^{**}_{H}(y_1,q_1,..., q_{N-1})- r&\leq& H_{L^{**}_H}(y_1, x_1,..., x_{N-1})  -H_{L^{**}_H}(y_2, x_1,..., x_{N-1})\\
& \leq& L^{**}_{H}(y_2,p_1,..., p_{N-1}) -L^{**}_{H}(y_1,p_1,..., p_{N-1}) +r.
\end{eqnarray*}
Since $L^{**}_{H}$ is Lipschitz,
\[-NR\|y_1-y_2\|- r\leq   H_{L^{**}_H}(y_1, x_1,..., x_{N-1})-H_{L^{**}_H}(y_2, x_1,..., x_{N-1}) \leq NR\|y_1-y_2\| +r.\]
  Since $r>0$ is arbitrary we obtain \[-NR\|y_1-y_2\|\leq   H_{L^{**}_H}(y_1, x_1,..., x_{N-1})-H_{L^{**}_H}(y_2, x_1,..., x_{N-1})  \leq NR\|y_1-y_2\|.\]
This together with ( \ref{lip}) prove that $H_{L^{**}_H}$ is Lipschitz continuous and that $Lip(H_{L^{**}_H})\leq NR.$
 \hfill $\square$\\

\textbf{Proof of Proposition \ref{prop1}.} 1) By Lemma \ref{ine1}, we have that $ H^1_{reg}:=H_{L^{**}_H}$ is a concave-convex Hamiltonian on $\R^d\times (\R^d)^{N-1}$ whose restriction to $\bar \Omega^N$ is $N$-sub-antisymmetric, hence  belong to $ {\cal H}^{^-}_N(\Omega)$.

2)\ To show that $H^2_{reg}$ is $N$-antisymmetric note that
\[
NH^2_{reg}({\bf x})=(N-1)H^1_{reg}({\bf x})-\sum_{i=1}^{N-1}H^1_{reg}(\sigma^{i}({\bf x}))=\sum_{i=1}^{N-1}\big[H^1_{reg}({\bf x})-H^1_{reg}(\sigma^{i}({\bf x}))\big]
\]
and each of the terms $H^1_{reg}({\bf x})-H^1_{reg}(\sigma^{i}({\bf x}))$ is easily seen to be $N$-antisymmetric.

Now $H^2_{reg}$ dominates $H^1_{reg}$ since
\[
N\big[H^2_{reg}({\bf x})-H^1_{reg}({\bf x})\big]=-H^1_{reg}({\bf x})-\sum_{i=1}^{N-1}H^1_{reg}(R^{i}({\bf x})) \geq 0,
\]
since $H^1_{reg}$ is $N$-sub-antisymmetric.

3)\ For $x\in \Omega$ and $p_1,..., p_{N-1} \in B_R$ we have
\begin{eqnarray*}
 L_{H^1_{reg}}(x, p_1,..., p_{N-1})&=& \sup\limits_{y_i\in \Omega}\Big \{\sum_{i=1}^{N-1}\langle p_i, y_i\rangle -H_{L^{**}_H}(x, y_1,..., y_{N-1})\Big\}\\
&=& \sup\limits_{y_i\in \Omega}\Big\{\sum_{i=1}^{N-1}\langle p_i, y_i\rangle -\sup\limits_{q_i\in B_R}\{\sum_{i=1}^{N-1}\langle q_i, y_i\rangle-L^{**}_H(x, q_1,..., q_{N-1}) \}\Big \}\\
&=& \sup\limits_{y_i\in \Omega}\inf \limits_{q_i\in B_R}\Big \{\sum_{i=1}^{N-1}\langle p_i, y_i\rangle -\sum_{i=1}^{N-1}\langle q_i, y_i\rangle+L^{**}_H(x, q_1,..., q_{N-1}) \}\Big \}\\
&\leq& \inf \limits_{q_i\in B_R} \sup\limits_{y_i\in \Omega}\Big \{\sum_{i=1}^{N-1}\langle p_i, y_i\rangle -\sum_{i=1}^{N-1}\langle q_i, y_i\rangle+L^{**}_H(x, q_1,..., q_{N-1}) \}\Big \}\\
&\leq&  \sup\limits_{y_i\in \Omega}\Big \{\sum_{i=1}^{N-1}\langle p_i, y_i\rangle -\sum_{i=1}^{N-1}\langle p_i, y_i\rangle+L^{**}_H(x, p_1,..., p_{N-1}) \}\Big \}\\
&=&L^{**}_H(x, p_1,..., p_{N-1})
\end{eqnarray*}
On the other hand by Lemma (\ref{ine2}) we have  $L^{**}_H \leq L_H $, and therefore $L_{H^1_{reg}} \leq L_H.$
It also follows from part 2) that $ L_{H^2_{reg}} \leq  L_{H^1_{reg}}.$
 This completes the proof of part 3).\\

Parts  4), 5) and 6) are  the subject of the preceding Lemmas. \hfill $\square$

\section{Proof of Theorem \ref{main}: Existence}

We first show that the minimization problem (${\rm DK}^2_{sym}$) is attained. Let $B_R$ be a ball such that $u_i(\bar \Omega) \subset B_R$ for all $i=1,..., N-1$.
Let $\{H_n\}$ be a sequence in ${\cal H}_N(\Omega)$ such that $L_{H_n}$ is a minimizing sequence for (${\rm DK}^2_{sym}$). Denoting $H_n^1:=(H_n)^1_{reg}$,
we get from Proposition \ref{prop1} that
$L_{H^1_n}\leq L_{H_n}$ on $\bar \Omega \times B^{N-1}_R$ and therefore $L_{H^1_n}$  %
is also minimizing for (${\rm DK}^2_{sym}$). It also follows from Proposition \ref{prop1} that
$L_{H^1_n}$ and $H^1_n$,
are uniformly  Lipschitz with  $Lip(H^1_n), Lip(L_{H^1_n}) \leq NR$. Moreover,
 \[
\hbox{$|H^1_n(x,y_1,..., y_{N-1})| \leq R\|x\|+R\sum_{i=1}^{N-1}\|y_i\| + 2NR^2$ for all $x$ and $(y_i)_{i=1}^{N-1}$ in $\R^d$,}
\]
and
 \[
\hbox{$ |L_{H^1_n}(x,p_1,..., p_{N-1})| \leq R\|x\|+R\sum_{i=1}^{N-1}\|p_i\| + (2N-1)R^2$  for all $x, p_1,..., p_{N-1}$ in $ \R^d.$}
 \]
  By Arzela-Ascoli's theorem, there exists two Lipschitz functions $\tilde H$ and  $\tilde L: \R^d \times \R^{d(N-1)} \to \R$ such that $H^1_n$
converges to $\tilde H$ and $L^1_n$ converges to $\tilde L$ uniformly on every compact set of $\R^N\times... \times \R^N.$  This implies that  $\tilde  H \in {\cal H}^{^-}_N(\Omega)$. Note that
\[L_{H^1_n}(x,p_1,..., p_{N-1}) +H^1_n(x, y_1,..., y_{N-1}) \geq  \sum_{i=1}^{N-1}\langle y_i,p_i\rangle, \]
for all $x,p_1,..., p_{N-1} \in \R^N $ and $y_1,..., y_{N-1} \in \bar \Omega,$ from which we have
 \[\tilde L(x,p_1,..., p_{N-1}) \geq \sum_{i=1}^{N-1}\langle y_i,p_i\rangle - \tilde H (x, y_1,..., y_{N-1}), \]
for all $x,p_1,..., p_{N-1} \in \R^N $ and $y_1,..., y_{N-1} \in \bar \Omega.$
This implies that $L_{\tilde H} \leq \tilde L$. Let $H_\infty^1=\tilde H^1_{reg}$ and $H_\infty^2=\tilde H^2_{reg}$ be the regularizations of $\tilde H $ defined in the
 previous section. Set $L^i_{\infty}=L_{H^i_{\infty}}$ for $i=1,2$. It  follows from Proposition \ref{prop1} that $L_{H^2_{\infty }}\leq L_{H^1_{\infty }} \leq L_{\tilde H}$ on
 $\bar \Omega \times B_R^{N-1},$ from which we have
\begin{eqnarray*}
{\rm DK}^2_{sym}&=&\int_{\Omega} L_{\tilde H}(x, u_1(x),..., u_{N-1}(x))\, d\mu\\
&=& \int_{\Omega} L^2_\infty(x, u_1(x),..., u_{N-1}(x))\, d\mu\\
&=&\int_{\Omega} L^1_{\infty}(x, u_1(x),..., u_{N-1}(x))\, d\mu.
\end{eqnarray*}
\hfill $\square$

For the rest of the proof, we shall need the following two technical lemmas. The first one relates $L_H^*$ to the standard Legendre transform of $H$ (extended beyond $\Omega^N$ to the whole of  $\R^{dN}$).

\begin{lemma} \label{wei} Let $H_\infty=H^1_\infty$ be the concave-convex Hamiltonian obtained above and $L_\infty=L^1_\infty$.  For each $x \in \bar \Omega$, define $
f_x: (\R^d)^{N-1}\to \R$ by
\[
f_x(y_1,..., y_{N-1}):=H_{\infty}(x, y_1,..., y_{N-1}).
\]
 We also define $\tilde f_x: (\R^d)^{N-1}\to \R \cup\{+\infty\}$ by
 \[
\hbox{$\tilde f_x (y_1,..., y_{N-1}):=f_x(y_1,..., y_{N-1})$ if $y_1,..., y_{N-1} \in \bar \Omega^{N-1}$ and $+\infty$ otherwise. }
\]
Let $(\tilde f_x)^{*}$ be the standard Fenchel dual of $\tilde f_x$ on $(\R^d)^{N-1}$ in such a way that  $(\tilde f_x)^{***}=(\tilde f_x)^{*}$ on $(\R^d)^{N-1}.$ We then have,
\begin{equation}\label{cheap}
\hbox{$f_x = (\tilde f_x)^{**}= \tilde f_x$ on $\bar \Omega^{N-1}$}
\end{equation}
and
\begin{eqnarray}\label{also}
L_{\infty}(x,p_1,..., p_{N-1})&=&
\sup\limits_{(z_i) \in \bar \Omega^{N-1}} \{\sum\limits_{i=1}^{N-1}\langle z_i,p_i
\rangle-(\tilde f_x)^{**}(z_1,..., z_{N-1})\}\nonumber\\
&=&\sup\limits_{(z_i) \in (\R^d)^{N-1}} \{\sum\limits_{i=1}^{N-1}\langle z_i,p_i
\rangle-(\tilde f_x)^{**}(z_1,..., z_{N-1})\}.
\end{eqnarray}
\end{lemma}
\textbf{Proof.}
Since $(\tilde f_x)^{**}$ is the largest convex function below $\tilde f_x$ we have  and $f_x \leq (\tilde f_x)^{**}\leq \tilde f_x,$ from which we obtain  $f_x = (\tilde f_x)^{**}= \tilde f_x$ on $\bar \Omega^{N-1}.$

For (\ref{also}), we first deduce from (\ref{cheap}) that
\begin{eqnarray*}
(\tilde f_x)^{*}(y_1,..., y_{N-1})&=& (\tilde f_x)^{***}(y_1,..., y_{N-1})\\
&= &\sup_{z \in \R^{d(N-1}} \{\sum\limits_{i=1}^{N-1}\langle z_i,y_i \rangle-(\tilde f_x)^{**}(z_1,..., z_{N-1})\}\\
&\geq& \sup_{z \in B^{N-1}_R} \{\sum\limits_{i=1}^{N-1}\langle z_i,y_i \rangle-(\tilde f_x)^{**}(z_1,..., z_{N-1})\}\\
&\geq&\sup_{z \in \Omega^{N-1}} \{\sum\limits_{i=1}^{N-1}\langle z_i,y_i \rangle-(\tilde f_x)^{**}(z_1,..., z_{N-1})\}\\
&=& \sup_{z \in \Omega^{N-1}} \{\sum\limits_{i=1}^{N-1}\langle z_i,y_i \rangle-f_x(z_1,..., z_{N-1})\}\\
&=& \sup_{z \in \Omega^{N-1}} \{\sum\limits_{i=1}^{N-1}\langle z_i,y_i \rangle- \tilde f_x(z_1,..., z_{N-1})\}\\
&=& (\tilde f_x)^{*}(y_1,..., y_{N-1}),
\end{eqnarray*}
from which we have the desired result.
 \hfill $\square$\\

 Fix now $H_\infty$ as above and let $H \in C({\bar \Omega}^N)$. For each  $\lambda >0$ and $r \in (-1,1) $, we  associate the following three functionals.
\begin{eqnarray*}
L_{r,\lambda}(x,p_1,..., p_{N-1})&:=& \sup_{(z_i) \in \bar \Omega^{N-1}}\left\{\sum_{i=1}^{N-1}\langle z_i,p_i
\rangle-(\tilde f_x)^{**}(z_1,..., z_{N-1}) -\frac{\lambda}{2}[\sum_{i=1}^{N-1}\|z_i\|^2- (N-1)\|x\|^2].\right. \\
&& \qquad \qquad \qquad \left. +rH(x, z_1,..., z_{N-1}) \right\} \\
L_{\lambda}(x,p_1,..., p_{N-1})&:=& \sup_{(z_i) \in  \R^{d(N-1)}}\left\{\sum_{i=1}^{N-1}\langle z_i,p_i
\rangle-(\tilde f_x)^{**}(z_1,..., z_{N-1}) -\frac{\lambda}{2}[\sum_{i=1}^{N-1}\|z_i\|^2- (N-1)\|x\|^2]  \right\}\\
L_{r}(x,p_1,..., p_{N-1})&:=& \sup_{(z_i) \in \bar \Omega^{N-1}}\left\{\sum_{i=1}^{N-1}\langle z_i,p_i
\rangle-H_{\infty}(x, z_1,..., z_{N-1})  +r
H(x, z_1,..., z_{N-1})\right\}.
\end{eqnarray*}

\begin{lemma} \label{limit} Let $H \in C({\bar \Omega}^N)$ be such that $H_\infty-rH \in {\cal H}^{^-}_N(\Omega)$ for all $r \in (-1,1). $ Then, the following  hold:
\begin{enumerate}
\item For every $(x, p_1,..., p_{N-1})\in \R^d\times \R^{d(N-1)}$, we have
\[
\hbox{$\lim\limits_{\lambda \to 0^+} L_{\lambda}(x,p_1,..., p_{N-1})=L_{\infty}(x,p_1,..., p_{N-1})$  and $\lim\limits_{\lambda \to 0^+} L_{r,\lambda}(x,p_1,..., p_{N-1})=L_{r}(x,p_1,..., p_{N-1}).$}
\]
\item  For all $x\in \R^d$, the function $(p_1,..., p_{N-1}) \to L_{\lambda}(x,p_1,..., p_{N-1})$ is differentiable.
\item For every $(x, p_1,..., p_{N-1})\in \R^d\times \R^{d(N-1)}$, we have
\[
\lim\limits_{r\rightarrow 0} \frac {L_{r, \lambda}(x,p_1,..., p_{N-1})-L_{\lambda}(x,p_1,..., p_{N-1})}{r}=H(\nabla_{2,...,N}
L_{\lambda}(x,p_1,..., p_{N-1}),x).
\]
\end{enumerate}
\end{lemma}
\textbf{Proof.} Yosida's regularization of convex functions and  Lemma \ref{wei}  yield that
\begin{eqnarray*}
\lim_{\lambda \to 0^+} L_{r,\lambda}(x,p_1,..., p_{N-1})&=&\sup_{(z_i) \in \bar \Omega^{N-1}}\{\sum_{i=1}^{N-1}\langle z_i,p_i
\rangle-(\tilde f_x)^{**}(z_1,..., z_{N-1}) -r
H(x, z_1,..., z_{N-1})\}\\
&=&\sup_{(z_i) \in \bar \Omega^{N-1}}\{\sum_{i=1}^{N-1}\langle z_i,p_i
\rangle-H_{\infty}(x, z_1,..., z_{N-1}) -r
H(x, z_1,..., z_{N-1})\}\\
&=&L_{r}(x,p_1,..., p_{N-1}).
\end{eqnarray*}
We also have
\[\lim_{\lambda \to 0}L_{\lambda}(x,p_1,..., p_{N-1})= \sup_{(z_i) \in \R^{d(N-1)}}\{\sum_{i=1}^{N-1}\langle z_i,p_i
\rangle-(\tilde f_x)^{**}(z_1,..., z_{N-1})\},\]
which, together with  Lemma \ref{wei}, yield that $\lim_{\lambda \to 0}L_{\lambda}(x,p_1,..., p_{N-1})=L_{\infty}(x,p_1,..., p_{N-1}).$

 (2) follows from the fact that the Yosida regularization of  convex functions are differentiable.

 (3) We let  $z_{_{(r, \lambda,i)}} \in \bar \Omega$ and $ z'_{_{(r, \lambda,i)}} \in \R^d$ be such that
\begin{eqnarray*}
L_{r,\lambda}(x,p_1,..., p_{N-1})&\leq& \sum_{i=1}^{N-1}\langle z_{_{(r, \lambda,i)}},p_i\rangle-(\tilde f_x)^{**}(z_{_{(r, \lambda,1)}}, ..., z_{_{(r, \lambda,N-1)}})-\frac{\lambda}{2}\sum_{i=1}^{N-1}\|z_{_{(r, \lambda,i)}}\|^2\\
&&+\lambda\frac{(N-1)\|x\|^2}{2}+r
H(x, z_{_{(r, \lambda,1)}},..., z_{_{(r, \lambda,N-1)}})+r^2,\\
L_{\lambda}(x,p_1,..., p_{N-1})&\leq& \sum_{i=1}^{N-1}\langle z'_{\lambda,i},p_i\rangle-
(\tilde f_x)^{**}(z'_{_{(r, \lambda,1)}}, ..., z_{_{(r, \lambda,N-1)}})-\frac{\lambda}{2}\sum_{i=1}^{N-1}\|z'_{_{(r, \lambda,i)}}\|^2+\lambda\frac{(N-1)\|x\|^2}{2}+r^2.
\end{eqnarray*}
Therefore,
\begin{eqnarray}\label{ine}
rH(x, z'_{_{(r, \lambda,1)}},..., z'_{_{(r, \lambda,N-1)}})-r^2&\leq& {L_{r,\lambda}(x,p_1,..., p_{N-1})-L_{\lambda}(x,p_1,..., p_{N-1})}\nonumber\\
&\leq&
rH(x, z_{_{(r, \lambda,1)}},..., z_{_{(r, \lambda,N-1)}})+r^2.
\end{eqnarray}
By the definition of  $L_\lambda,$ we have  $\sup_{r \in [-1,1]}\|z'_{r, \lambda,i}\| < \infty.$ Suppose now that, up to a subsequence, $z_{r,\lambda,i}\rightarrow z_{i} \in \bar \Omega$ and $z'_{r, \lambda,i}\rightarrow z_{\lambda,i}'$ as $r \to 0.$  This together with the definition of
$L_{r,\lambda}$ and $L_{\lambda}$ imply that
\begin{eqnarray*}\label{equality}
L_{\lambda}(x,p_1,..., p_{N-1})&=& \sum_{i=1}^{N-1}\langle z_{_{(\lambda,i)}},p_i\rangle-(\tilde f_x)^{**}(z_{_{(\lambda,1)}},..., z_{_{(\lambda,N-1)}})-\frac{\lambda}{2}\sum_{i=1}^{N-1}\|z_{_{(\lambda,i)}}\|^2+\lambda (N-1)\frac{\|x\|^2}{2}\\
&=&\sum_{i=1}^{N-1}\langle z'_{\lambda,i},p_i\rangle -(\tilde f_x)^{**}(z_{_{(\lambda,1)}},..., z'_{_{(\lambda,N-1)}})-\frac{\lambda}{2}\sum_{i=1}^{N-1}\|z'_{_{(\lambda,i)}}\|^2+\lambda (N-1)\frac{\|x\|^2}{2},
\end{eqnarray*}
from which we obtain that
\begin{equation}\label{omeg}z_{\lambda,i}=z'_{\lambda,i}=\nabla_i L_{\lambda}(x,p_1,..., p_{N-1}) \in \bar \Omega, \qquad \quad i=2,..., N.\end{equation}
It then follows
from (\ref{ine}) that
\begin{eqnarray*}
\lim_{r\rightarrow 0} \frac {L_{r, \lambda}(x,p_1,..., p_{N-1})-L_{\lambda}(x,p_1,..., p_{N-1})}{r}=H\big(\nabla_{2,..., N}
L_{\lambda}(x,p_1,..., p_{N-1}),x\big).
\end{eqnarray*}
\hfill$\square$ \\
\textbf{End of the proof of Theorem \ref{main} (Existence):}  For each $\lambda >0, $ $x \in \bar \Omega$ and $p \in \R^N$, we define
\[
\bar S_{\lambda, i} (x,p_1,...,p_{N-1})=\nabla_i L_{\lambda}(x,p_1,..., p_{N-1})\quad i=2,..., N
\]
We have that $\bar S_{\lambda, i} (x,p_1,...,p_{N-1}) \to \bar S_{0, i} (x,p_1,...,p_{N-1})$ where $\bar S_{0, i} (x,p_1,...,p_{N-1})$ is the unique
element with minimal norm in $\partial_{i} L_{\infty}(x,p_1,..., p_{N-1})$.\\
  Set $S_{\lambda,i}(x)= \bar S_{\lambda,i}(x, u_1(x),...,u_{N-1}(x))$ and
 $S_i(x)= S_{0,i}(x, u_1(x),..., u_{N-1}(x)).$ For each $r>0,$ $\lambda \in [0,1]$ and $x \in \bar \Omega,$ define
\[\eta_r(\lambda, x)= \frac{L_{r,\lambda}(x, u_1(x),..., u_{N-1}x) -L_{\lambda}(x, u_1(x),..., u_{N-1}(x))}{r}.\]
Note that the function $r \to L_{r,\lambda}(x, u_1(x),..., u_{N-1}(x))$ is a convex function because it is supremum of a family of linear functions. Thus, for fixed $(x, \lambda)\in \Omega \times [0,1]$,  the function $r \to \eta_r(\lambda, x)$ is non-decreasing. Setting $\eta_0(\lambda, x)$ to be $H(x, S_{\lambda,1}(x),..., S_{\lambda,N-1}(x))$ for $\lambda >0$ and $\eta_0(0, x)=H(x, S_{1}(x),..., S_{N-1}(x))$, we have that both  functions $\lambda \to \eta_r(\lambda, x)$ and $\lambda \to \eta_0(\lambda, x)$ are continuous. It follows from Dini's Theorem, that for a fixed $x,$ $\eta_r(\lambda, x)$ converges uniformly to $\eta_0(\lambda, x)$ as $r \to 0$ with respect to $\lambda \in [0,1].$ Note also that thanks to (\ref{omeg}) we have that  $S_{\lambda, i}, S_i:\bar \Omega \to \bar \Omega$ and  for all $x\in \Omega$.
\begin{equation}\label{key}
(S_1x,..., S_{N-1}x) \in \partial_{2,...,N} L_{\infty}(x, u_1(x),..., u_{N-1}(x)).
\end{equation}
We now show that
\begin{equation}\label{admiss}
\hbox{$\int_{\Omega} H(x, S_1x,..., S_{N-1}x) \, d\mu=0$ for all $H \in C({\bar \Omega}^N)$ with $H_\infty-rH \in {\cal H}^{^-}_N(\Omega)$,  $r \in (-1,1)$.}
\end{equation}
 Indeed, since $|H(x, S_1x,..., S_{N-1}x) | \leq \|H\|_{L^\infty({\bar \Omega}^N)},$ we get from Lebesgue's dominated convergence Theorem,
\[\lim_{\lambda \to 0} \int_{\Omega} H(x, S_{\lambda,1}(x),..., S_{\lambda,N-1}(x)) \, d\mu (x) =\int_{\Omega}  H(x, S_1x,..., S_{N-1}x) \, d\mu (x).\]
From (\ref{ine}) we have
\[\Big |\frac {L_{r, \lambda}(x,p_1,..., p_{N-1})-L_{\lambda}(x,p_1,..., p_{N-1})}{r} \Big| \leq \|H\|_{L^\infty({\bar \Omega}^N)}+|r|, \]
from which follows that
\begin{eqnarray*}
 \int_{\Omega} H(x, S_1x,..., S_{N-1}x) \, d\mu (x) &=&\int_{\Omega} \lim_{\lambda \to 0}  \lim_{r\rightarrow 0^+} \frac {L_{r, \lambda}(x,u_1(x),..., u_{N-1}(x))-L_{\lambda}(x,u_1(x),..., u_{N-1}(x))}{r}\, d\mu\\
&=&\int_{\Omega} \lim_{\lambda \to 0}  \lim_{r\rightarrow 0^+} \eta_r(\lambda, x)\, d\mu\\
&=&\int_{\Omega}  \lim_{r\rightarrow 0^+} \lim_{\lambda \to 0}  \eta_r(\lambda, x)\, d\mu \qquad  \text{(due to the uniform convergence) }\\
&=&\int_{\Omega}  \lim_{r\rightarrow 0^+}   \eta_r(0, x)\, d\mu \\
&=&\lim_{r\rightarrow 0^+} \int_{\Omega}    \eta_r(0, x)\, d\mu  \qquad \text{ (due to the monotone convergence theorem) } \\
&=&  \lim_{r\rightarrow 0^+} \int_{\Omega} \frac {L_{r}(x,u_1(x),..., u_{N-1}(x))-L_{\infty}(x,u_1(x),..., u_{N-1}(x))}{r}\, d\mu\\
&\geq & 0,\qquad \quad \text{(in view of the optimality of $H_\infty $ compared to $H_\infty -r H$). }
\end{eqnarray*}
In other words, we have $ \int_{\Omega} H(x, S_1x,..., S_{N-1}x) \, d\mu \geq 0$. By the same argument considering $r \to 0^-$, one has
 $ \int_{\Omega} H(x, S_1x,..., S_{N-1}x) \, d\mu \leq 0$ and therefore  the latter is indeed zero as desired.

 Note now that (\ref{admiss}) yields that
 \begin{equation}\label{admi}
 \int_{\Omega} H_{\infty}(x, S_1x,..., S_{N-1}x) \, d\mu= 0,
 \end{equation}
and
\begin{equation}
\hbox{$ \int_{\Omega} H(x, S_1x,..., S_{N-1}x) \, d\mu = 0$ for all $H\in {\cal H}_N(\Omega)$.}
 \end{equation}
 It follows from Lemma \ref{magic} that $S$ is measure preserving, that $S_i=S_1^{i}$ and that $S_1^{N}=I$. We shall now write $S$ for $S_1$.

In order to show that ${\rm DK}^2_{\rm sym}={\rm MK}_{\rm sym}$, we note that clearly ${\rm MK}_{\rm sym} \leq {\rm DK}^2_{\rm sym}$. For the reverse inequality, we use the fact that $(Sx,..., S^{N-1}x) \in \partial_{2,...,N} L_{\infty}(x, u_1(x),..., u_{N-1}(x))$ together with  $(\tilde f_x)^{**}$ being  the Fenchel dual of $L$ with respect to the last $N-1$ variables and Lemma \ref{wei} to obtain that
\begin{equation} \label{inter1}
(u_1(x),..., u_{N-1}(x))\in \partial (\tilde f_x)^{**} (Sx, ..., S^{N-1}x).
\end{equation}
Since
${\rm meas}(\partial \Omega)=0,$ the set $\displaystyle{\cup_{i=1}^{N-1}
S^{-i}(\partial \Omega)}$ is  negligible and for each $x \in \Omega\setminus \cup_{i=1}^{N-1}
S^{-i}(\partial \Omega)$,
one has
\[
\partial (\tilde f_x)^{**} (Sx, ..., S^{N-1}x)=\partial_{2,...,N} H_\infty (x, Sx, ..., S^{N-1}x).
\]
 It follows that
\begin{equation}\label{inter2}
(u_1(x),..., u_{N-1}(x))\in \partial_{2,...,N} H_\infty (x, Sx, ..., S^{N-1}x)\qquad a.e. \quad x \in \Omega.
\end{equation}
We finally get that
\begin{eqnarray*}
{\rm DK}^2_{\rm sym}&=& \int_{\Omega} L_{\infty}(x, u_1(x),..., u_{N-1}(x)) \, d\mu (x)\\
 &=&\int_{\Omega} L_{\infty}(x, u_1(x),..., u_{N-1}(x)) \, d\mu (x)+ \int_{\Omega}H_{\infty}(x, Sx, ..., S^{N-1}x) \, d\mu (x)\\&=&\int_{\Omega} L_{\infty}(x, u_1(x),..., u_{N-1}(x)) \, d\mu (x)+ \int_{\Omega}(\tilde f_x)^{**}(Sx, ..., S^{N-1}x) \, d\mu (x)
\\&=& \int_{\Omega} \sum_{i=1}^{N-1}\langle u_i(x), S^i (x) \rangle \, d\mu (x)\leq {\rm MK}_{\rm sym}.
\end{eqnarray*}
If now $u_i \in W^{1,1}_{loc} (\Omega)$ for $i=1,2,...,N-1,$ or if $S$ is a.e. differentiable, then by Theorem \ref{dar} of the Appendix, there exists a full measure subset $\Omega_0$ of $\Omega$ that
$\nabla_{2,...,N} H_{\infty}(x, Sx,..., S^{N-1}x)$ exists for all $x \in  \Omega_0.$ It follows that
\[\big (u_1(x),..., u_{N-1}(x)\big )=\nabla_{2,...,N} H_{\infty}(x, Sx, ..., S^{N-1}x) \quad \text{ for all } x \in \Omega_0.\]

\section{Proof of Theorem \ref{main}: Uniqueness}
We now deal with part (5) of Theorem \ref{main}. $H_\infty$ will denote an optimal concave-convex $N$-sub-antisymmetric associated to the vector fields $u_1,..., u_{N-1}$ obtained via the above variational procedure.
\begin{lemma} \label{H8} Assume that the vector fields $u_1,..., u_{N-1}$ from $\Omega$ to $\R^d$ are such that
\[
\hbox{$\big (u_1(x),..., u_{N-1}(x)\big )\in \partial_{2,...,N} H_{1}(x, Sx, ..., S^{N-1}x)$ \, a.e. $x\in \Omega$},
\]
for some concave-convex $N$-sub-antisymmetric Hamiltonian $H_1$ and some $N$-involution $S$, then  $(H_1, S)$ is an ``extremal pair", meaning that the infimum $({\rm DK}^2_{\rm sym})$ is attained at $H_1$ and the supremum $({\rm MK}_{cyc})$ is attained at $S$. Moreover, we have
\[
\hbox{$\big (u_1(x),..., u_{N-1}(x)\big )\in \partial_{2,...,N} H_{\infty}(x, Sx, ..., S^{N-1}x)$ \, a.e. $x\in \Omega$},
\]
  where $H_\infty$ is the optimal Hamiltonian constructed in Theorem \ref{main}.

\end{lemma}
\textbf{Proof.} Let $L$ be the Fenchel-Legendre dual of $H_1$ with respect to the last $N-1$ variable. We have that $L_{H_1}\leq L$ on $(\R^d)^{N-1} \times \Omega$. It follows that
\begin{eqnarray*}\sum_{i=1}^{N-1}\langle u_i(x), S^i (x) \rangle  &\leq &L_{H_1}(x,u_1(x),..., u_{N-1}(x))+ H_1(x, S x,..., S^{N-1}x)\\
 &\leq &L(x,u_1(x),..., u_{N-1}(x))+ H_1(x, S x,..., S^{N-1}x)\\
 &=&
\sum_{i=1}^{N-1}\langle u_i(x), S^i (x) \rangle,
\end{eqnarray*}
from which we deduce that
\[\sum_{i=1}^{N-1}\langle u_i(x), S^i (x) \rangle  =L_{H_1}(x,u_1(x),..., u_{N-1}(x))+ H_1(x, S x,..., S^{N-1}x),\]
and \[\int_\Omega \sum_{i=1}^{N-1}\langle u_i(x), S^i (x) \rangle\, d\mu= \int_\Omega L_{H_1}(x,u_1(x),..., u_{N-1}(x)) \, d\mu+\int_\Omega H_1(x, S x,..., S^{N-1}x)\, d\mu.\]
Use now the optimality of $H_1 $ compared to $H_1-r H_1$ for $-1<r<1$ (Indeed, the above equality will be an inequality when $H_1$ is replaced by $H_1-r H_1$ for $r \not =0$) and the same argument  as in the proof of the existence part in Theorem \ref{main} for $H_\infty$ to obtain that
$\int_\Omega H_1(x, S x,..., S^{N-1}x)\, d\mu=0$.
On the other hand, we have
\[\int_\Omega \sum_{i=1}^{N-1}\langle u_i(x), S^i (x) \rangle \, d\mu \leq {\rm MK}_{cyc}={\rm DK}^2_{\rm sym} \leq \int_\Omega L_{H_1}(x,u_1(x),..., u_{N-1}(x)) \, d\mu,\]
which yields
\[\int_\Omega \sum_{i=1}^{N-1}\langle u_i(x), S^i (x) \rangle  \, d\mu = {\rm MK}_{cyc}={\rm DK}^2_{\rm sym} = \int_\Omega L_{H_1}(x,u_1(x),..., u_{N-1}(x)) \, d\mu.
\]
Now we can  show that $u_i(x)\in \partial_{i+1}H_\infty(x, S x,..., S^{N-1}x)$ a.e.  In fact,
\begin{eqnarray*}\int_\Omega \sum_{i=1}^{N-1}\langle u_i(x), S^i (x) \rangle \, d\mu &=&\int_\Omega L_{H_1}(x,u_1(x),..., u_{N-1}(x)) \, d\mu\\
&=&{\rm DK}^2_{\rm sym}=
\int_\Omega L_{H_\infty}(x,u_1(x),..., u_{N-1}(x)) \, d\mu \\&\geq& \int_\Omega L_\infty(x,u_1(x),..., u_{N-1}(x)) \, d\mu +\int_\Omega H_\infty(x, S x,..., S^{N-1}x) \, d\mu\\
&\geq &\int_\Omega \sum_{i=1}^{N-1}\langle u_i(x), S^i (x) \rangle \, d\mu,
\end{eqnarray*}
which implies that
\[
\hbox{$\sum_{i=1}^{N-1}\langle u_i(x), S^i (x) \rangle =L_{H_\infty}(x,u_1(x),..., u_{N-1}(x)) +H_\infty(x, S x,..., S^{N-1}x)$ a.e. on $\Omega$, }
\]
and hence the desired result. \hfill $\Box$

\begin{lemma}\label{H9} Suppose $S$ is a measure preserving N-involution and $u_i(x) = \nabla_{i+1}H_\infty(x, S x,..., S^{N-1})$ a.e. for $i=1,..., N-1$. Then
\[\nabla_1 H_\infty(x, S x,..., S^{N-1}x)=-\sum_{i=1}^{N-1}u_i(S^{N-i}x) \qquad \text{ a.e. } x \in \Omega.\]
 \end{lemma}
\textbf{Proof.}
Let $u \in \R^d$ and let $|t|$ be small.  Note that
\[\int_\Omega \sum_{i=1}^{N} H_\infty \big (\sigma^{N+1-i} (x, Sx, ..., S^{N-1}x)\big)\, d\mu = N\int_\Omega H_\infty \big( x, Sx, ..., S^{N-1}x\big) \, d\mu =0.\]
Since  $\sum_{i=1}^{N} H_\infty(\sigma^{N+1-i}\big (x, Sx, ..., S^{N-1}(x)\big)) \leq 0$, it follows that
\[\sum_{i=1}^{N} H_\infty\big (\sigma^{N+1-i} (x, Sx, ..., S^{N-1}x)\big) =0 \quad \text{ a.e. } x \in \Omega. \]
Note that $H_\infty$ is $N$-sub-antisymmetric and therefore
\begin{eqnarray*}
 \sum_{i=1}^{N} H_\infty\big (\sigma^{N+1-i} (x+tu, Sx, ..., S^{N-1}x)\big) \leq 0 = \sum_{i=1}^{N} H_\infty\big (\sigma^{N+1-i} (x, Sx, ..., S^{N-1}x)\big).
\end{eqnarray*}
Assuming $x$ is a point where  $\nabla_i H_\infty\big (\sigma^{N+1-i} (x+tu, Sx, ..., S^{N-1}x)\big)$ exists for all $i=1,...,N,$ then
\[
  \sum_{i=1}^{N} \nabla_i H_\infty \big (\sigma^{N+1-i}(x, Sx, ..., S^{N-1}x)\big)=0.
\]
Since  $u_i(x)=\nabla_{i+1} H_\infty(x, Sx,.., S^{N-1}x)$ and $S^N=I$ a.e., we have for $i=2,3,...,N$,
\[u_{i-1}(S^{N+1-i} x)=\nabla_{i} H_\infty\big (\sigma^{N+1-i} (x, Sx, ..., S^{N-1}x)\big).
\] Therefore,
\[
  \sum_{i=1}^{N-1}u_i(S^{N-i} x)+\nabla_1 H_\infty(x, S x,..., S^{N-1}x) =0.
  \]
\begin{proposition} \label{critical} Let $u_1,..., u_{N-1}$ be vector fields in $W^{1,1}_{loc} (\Omega)$ such that for any two families of points $x_1,..., x_N$ and $y_1,..., y_N$ in $\Omega$,
 the function
\[x\to\sum_{i=1}^{N-1}\langle u_i(x), y_i-x_i\rangle +\sum_{i=1}^{N-1}\langle u_i(y_{N-i})-u_i(x_{N-i}), x\rangle\]
has no critical point unless when $x_1=y_1.$ Then, there is a unique measure preserving $N$-involution $S$ on $\Omega$ that satisfies $(\ref{rep})$ for some concave-convex $N$-sub-antisymmetric Hamiltonian $H$.
\end{proposition}
\textbf{Proof.} Suppose $S_1$, $S_2$ are two measure preserving $N$-involutions on $\Omega$ and $H_1$ and $H_2$ are two concave-convex $N$-sub-antisymmetric Hamiltonian on $\Omega \times \Omega^{N-1}$ such that for $j=1,2$, we have
\begin{equation}
 u_i(x)=\nabla_i H_j (x, S_j^1x,..., S_j^{N-1}) \quad i=1,..., N-1.
\end{equation}
Note first that Lemma \ref{H8} gives that
\begin{equation}
 u_i(x)=\nabla_i H_\infty (x, S_j^1x,..., S_j^{N-1}).
\end{equation}
From Lemma \ref{H9}, we have that
\[
-\sum_{i=1}^{N-1}u_i(S_j^{N-i} x)=\nabla_1 H_\infty (x, S_j^1x,..., S_j^{N-1}).
\]
Note that the function $x\to L_\infty(x,u_1,..., u_{N-1}(x))$ is locally Lipschitz and therefore is differentiable on a subset $\Omega_0$ of full measure.
We now show that $S_1=S_2$ on $\Omega_0$.

Indeed, for any $x\in \Omega_0$, $h=0$ is a minimum for the function
\begin{eqnarray*}
h\to L_\infty(x+h, u_1(x+h),..., u_{N-1}(x+h)) +H_\infty( x+h,S_j^1x,..., S_j^{N-1}x )-\sum_{i=1}^{N-1}\langle u_i(x+h), S^i_j(x)\rangle.
\end{eqnarray*}
This implies that
\begin{eqnarray*}
\nabla_{1}H_\infty ( x,S_1^1x,..., S_1^{N-1}x) -\sum_{i=1}^{N-1}\langle \nabla u_i(x), S^i_1(x)\rangle&=&-\frac{d}{dh}L_\infty(x+h,u_1(x+h),..., u_{N-1}(x+h))_{h=0}\\
&=&\nabla_{1}H_\infty ( x,S_2^1x,..., S_2^{N-1}x) -\sum_{i=1}^{N-1}\langle \nabla u_i(x), S^i_2(x)\rangle.
\end{eqnarray*}
This yields that
\begin{eqnarray*}
\sum_{i=1}^{N-1}\langle \nabla u_i(x), S^i_2(x)-S^i_1(x)\rangle&=& \nabla_{1}H_\infty ( x,S_2^1x,..., S_2^{N-1}x)-\nabla_{1}H_\infty ( x,S_1^1x,..., S_1^{N-1}x)\\&=&
\sum_{i=1}^{N-1}\Big (u_i(S^{N-i}_1(x))-u_i(S^{N-i}_2(x)\Big ).
\end{eqnarray*}
The hypothesis then implies that $S_1(x)=S_2(x)$, and $S$ is therefore unique.\\

In order to find examples of families of vector fields satisfying the above sufficient condition for uniqueness, we look again at $N$-monotone vector fields.   For that we introduce the following notion.

\begin{definition} Say that a family of vector fields  $(u_1, u_2,..., u_{N-1})$ on $\Omega$  is {\it strictly jointly $N$-monotone} if
\begin{equation}\label{uni0} \sum_{i=1}^{N}\sum_{l=1}^{N-1}\langle u_l(x_i), x_i-x_{i+l} \rangle > 0,\end{equation}
for every cycle $x_1,..., x_{2N-1}$ of points in $\Omega$
such that $x_{N+l}=x_l$ for $1 \leq l \leq N-1$, and $x_1\neq x_2$.
\end{definition}
Note that for $N=2$, this property means that the vector field $u_1$ is {\it  strictly $2$-monotone}, that is,
\begin{eqnarray}\label{uu}
\hbox{$ \langle u_1(y)-u_1(x), y-x\rangle > 0$ for all $y, x \in \Omega$ with $x\not=y.$}
\end{eqnarray}
In this case, it is easy to see that if $u_1$ is  differentiable, then strict monotonicity implies the sufficient condition for uniqueness mentioned in Proposition \ref{critical}. Indeed, let $u \in \R^d$ and $x \in \Omega.$ By taking $y=x+tu$  in (\ref{uu}) and letting $t \to 0^+$ we obtain
$\langle \nabla u_1(x) u, u \rangle \geq 0.$

Assume now that the function
$
x\to\langle u_1(x), y_1-x_1\rangle +\langle u_1(y_1)-u_1(x_1), x\rangle
$
has a  critical point and that $y_1\not=x_1$. It follows that
\[
\langle \nabla u_1(x) (y_1-x_1) ,y_1-x_1\rangle +\langle u_1(y_1)-u_1(x_1), y_1-x_1\rangle=0.
\]
Since the first term is non-negative and the second one is strictly positive, this leads to a contradiction.

One can however, establish directly the following uniqueness result for strictly jointly $N$-monotone families for $N\geq 2$, even  without the differentiability assumption on $u_1, ..., u_{N-1}.$ This is because we already know from the result of Galichon-Ghoussoub \cite{GG} mentioned in the introduction that  $S_1(x)=x$ is one of the possible $N$-involution measure preserving maps  in the representation of $(u_1, ..., u_{N-1})$.

\begin{proposition} Assume $u_1, ..., u_{N-1}$ is a strictly jointly $N$-monotone family of bounded vector fields on $\Omega$. Then, $S=I$ is the only  measure preserving $N$-involution $S$ on $\Omega$ that satisfies $(\ref{rep})$ for some concave-convex $N$-sub-antisymmetric Hamiltonian $H$.
\end{proposition}
\textbf{Proof.}  Assume $S$ is  another measure preserving $N$-involution in the decomposition. Let $x_i=S^ix$ for $i=1,2,..., N$ and note that $x_N=x.$
It follows from (\ref{uni0}) that
\[\sum_{i=0}^{N-1}\sum_{l=1}^{N-1}\langle u_l(S^ix), S^ix-S^{i+l}x \rangle \geq 0.
\]
Integrating the above expression over $\Omega$ implies that
\begin{eqnarray*}
 0 &\leq&\int_\Omega \sum_{i=0}^{N-1}\sum_{l=1}^{N-1}  \langle u_l(S^ix), S^ix-S^{i+l}x \rangle \, d\mu\\
 &=&\sum_{i=0}^{N-1}\sum_{l=1}^{N-1} \int_\Omega \langle u_l(S^ix), S^ix \rangle \, d\mu-\sum_{i=0}^{N-1}\sum_{l=1}^{N-1} \int_\Omega \langle u_l(S^ix), S^{i+l}x \rangle \, d\mu\\
  &=&\sum_{i=0}^{N-1}\sum_{l=1}^{N-1} \int_\Omega \langle u_l(x), x \rangle \, d\mu-\sum_{i=0}^{N-1}\sum_{l=1}^{N-1} \int_\Omega \langle u_l(x), S^{l}x \rangle \, d\mu\\
   &=&N\sum_{l=1}^{N-1} \int_\Omega \langle u_l(x), x \rangle \, d\mu-N\sum_{l=1}^{N-1} \int_\Omega \langle u_l(x), S^{l}x \rangle \, d\mu\\
   &=&N\int_\Omega L_\infty(x,u_1(x),..., u_{N-1}(x)) \, d\mu -N\int_\Omega L_\infty(x,u_1(x),..., u_{N-1}(x)) \, d\mu \\
   &=&0.
\end{eqnarray*}
The latter identity is because both terms correspond to the optimal value (${\rm MK}_{cyc}$).   Since the integrand in the first line of the above expression is nonnegative we obtain
 \[\sum_{i=0}^{N-1}\sum_{l=1}^{N-1}\langle u_l(S^ix), S^ix-S^{i+l}x \rangle = 0, \qquad \text{ a.e. } x\in \Omega,
 \]
 and therefore $Sx=x.$ \hfill $\square$

\section{Proof of Theorem \ref{main.2}}

The question here is what happens when some of the vector fields  $u_i$ are identically zero.
Let us illustrate the situation by assuming that just one of them, say  $u_{N-1}\equiv 0.$
In this case, there are two scenarios:

{\bf (I)} \,\, One can begin  with $N-2$ vectors $u_1,.., u_{N-2},$ and obtain a sub $(N-1)-$antisymmetric Hamiltonian $H$ and an $(N-1)-$involution $S$ such that
$u_i(x) \in \partial H_{i+1}(x, Sx,..., S^{N-2}x).$\\

{\bf (II)} One can proceed as above, while considering $u_{N-1}\equiv 0$ as a vector field like the others. Note that in the proof of the main theorem we never assumed $u_{N-1} \not =0,$ except on line  (\ref{inter2}) and the preceding  paragraph. However, it is easily seen that
by assuming $u_{N-1}\equiv 0,$ one still gets
\begin{equation*}\label{inter3}
(u_1(x),..., u_{N-2}(x))\in \partial_{2,...,N-2} H_\infty (x, Sx, ..., S^{N-1}x)\qquad a.e. \quad x \in \Omega.
\end{equation*}
and the dependence of the Hamiltonian $H_\infty$ with respect to the $N$-th variable seems to be redundant. In this case  $H_\infty$ can be chosen to be an $N-$antisymmetric Hamiltonian, which depends on only $N-1$ variables. This is because  $H^1_{reg}(x_1,..., x_N)= H_{L^{**}_H}(x_1,,..., x_{N-1}, x_N)$ can be replaced by
\begin{eqnarray}\label{wha3} H^0_{reg}(x_1,..., x_{N-1}, x_N):=F_0(x_1,..., x_{N-1}), \end{eqnarray} where
\[F_0(x_1,,..., x_{N-1})=\sup_{p_2,..., p_{N-1} \in B_R}\Big \{ \sum_{i=2}^{N-1}\langle p_i, x_i\rangle- L^{**}_H(x_1,p_2,..., p_{N-1},0 ) \Big \}.\]
Indeed, it follows from (\ref{wha}) that  for all $x_i \in \Omega $ and $p_j^i \in \R^d$ the following inequality holds
\begin{eqnarray*}
 \sum_{i=1}^{N}\Big \{ \sum_{j=1, j\not =i}^{N}\langle p_j^i, x_j\rangle- L^{**}_H(R^{i-1}(p^i_1,..., p_{i-1}^i,x_i, p_{i+1}^i,..., p^i_{N}) ) \Big \}\leq 0.
\end{eqnarray*}
  In the above expression,  set   $p^i_{i-1}=p^1_N=0$ for $i>1$.  By taking sup over all non-zero  $p_i^j \in B_R$
we obtain
\begin{equation}\label{wha2}\sum_{i=1}^{N} H^0_{reg}\big (\sigma^{i-1}(x_1,...,x_N)\big ) \leq 0.\end{equation}
This proves that $H^0_{reg}$ is  $N-$sub-antisymmetric.
  By defining
\begin{equation*}
H^2_{reg}({\bf x})=\frac{(N-1)H^0_{reg}({\bf x})-\sum_{i=1}^{N-1}H^0_{reg}(\sigma^{i}({\bf x}))}{N},
\end{equation*}
and using a similar argument as in the proof of Proposition \ref{prop1}, one can  also obtain that   $L_{ H^2_{reg}}\leq L_{ H^0_{reg}} \leq L_{ H}$ on $\bar \Omega \times (B_R)^{N-1}$.
 This  together with ( \ref{wha3}) and ( \ref{wha2}) imply  that  the Hamiltonian $H_\infty$ obtained variationally in Theorem \ref{main} can be chosen to be independent with respect to the last variable.
\\

Similarly, one can show that if more than one vector fields is zero, then the dependence of $H_\infty$ on the corresponding variables can be dropped.

Suppose now that $u_2=...=u_{N-1}=0$. In this case $H_\infty$ is just a function of two variables, i.e. $H(x_1,x_2,...,x_N)=F(x_1,x_2)$
for some Lipschitz function $F$, which is concave with respect to the first variable and convex with respect to the second one. Therefore $u_1(x) \in \partial_2 F(x, Sx)$ for some
measure preserving $N-$involution.  In this case,
the sub-N-antisymmetry of $H_\infty $ translates into
\begin{eqnarray*}
\hbox{$ \sum_{i=1}^N F(x_{i+1}, x_i) \leq 0$ for all $x_1,...x_N \in \Omega $ with $x_1=x_{N+1}.$   }
\end{eqnarray*}

 \section{Appendix}

\begin{theorem} \label{dar} Consider bounded vector fields $(u_i)_{i=1}^{N-1}$ on $\Omega$ such that for $i=1,2,..., N-1$,
\begin{equation}
\hbox{$u_i (x)\in \partial_{i+1} H(x, Sx,..., S^{N-1}x)$ a.e. $\Omega$,}
\end{equation}
 where $S: \bar \Omega \to \bar \Omega$ is a measure preserving $N$-involution, and $H: \R^d \times (\R^d)^{N-1}$ is a Lipschitz function satisfying  the following properties:
\begin{enumerate}
\item  $H(\, .\,, X)$ is concave for every $X \in (\R^d)^{N-1},$  and $H(x,\, . \, )$ is convex for all $x \in \R^d.$
\item $H$ is $N$-sub-antisymmetric on $(\bar \Omega)^N.$
\item  $\int_{\Omega} H(x, Sx,..., S^{N-1}x) \, d\mu=0$.
 \end{enumerate}
If either $S \in W_{loc}^{1,1}(\Omega))$ or $u_i\in W_{loc}^{1,1}(\Omega)$
for $i=1,2,..., N-1$,
 then  there exists a full measure subset $\Omega_0$ of $\Omega$ such that
$\nabla_i H(x, Sx,..., S^{N-1}x)$ exists for all $x \in \Omega_0.$
\end{theorem}
First recall the following standard lemma.

\begin{lemma}\label{rock2} Let $f:\R^n \to (-\infty, +\infty]$ be a proper convex function and let  $x$ be a point where $f$ is finite. The following statements hold:
\begin{enumerate}
\item  For each $v \in \R^n,$ the difference quotient in the definition of  $Df(x)v$ is a  non-decreasing function of $\lambda >0,$ so that  $Df(x)v$  exists and
\begin{equation}\label{subd} Df(x)v= \inf_{\lambda>0} \frac{f (x+\lambda v)-f(x)}{\lambda}.\end{equation}
\item the function $v \to Df(x)v$ is a positively homogeneous convex function of $v$ with
\[Df(x)u+Df(x)(-v)\geq 0 \qquad \forall v \in \R^n.\]
\end{enumerate}
\end{lemma}

\begin{lemma} \label{seq0}For each $v \in \R^d,$ we have
\[\int_\Omega D_1 H (x, Sx, S^2x,...,S^{N-1}x)(v)  \, d\mu+\int_Q D_1 H (x, Sx, S^2x,...,S^{N-1}x)(-v)  \, d\mu = 0.\]
\end{lemma}
\textbf{Proof.}
Let  $t>0$ and define
\begin{eqnarray*}
I^1(x,v,t)&=& H (x, S(x+tv), S^2(x+tv),...,S^{N-1}(x+tv)),\\
I^2(x,v,t)&=& H (x+tv, Sx, S^2x,...,S^{N-1}x).
 \end{eqnarray*}
 Let  $g \in C_{c}^{\infty}(\Omega)$ be a non-negative  function. By a simple change of variables, we have for $t>0$ small enough,
\begin{eqnarray}\label{fin}
\int_{\Omega}\frac{I^1(x,v,t)g(x)+I^1(x,-v,t)g(x)-2I^1(x,0,0)g(x)}{t}\, d\mu= \nonumber\\
\int_{\Omega}\frac{I^2(x,-v,t)g(x-t v)+I^2(x,v,t)g(x+tv)-2I^1(x,0,0)g(x)}{t}\, d\mu.
\end{eqnarray}
The limit of the right hand side of the above expression exists as $t \to 0^+$ and
\begin{eqnarray}\label{fin1}
\lim_{t \to 0^+}\int_{\Omega}\frac{I^2(x,-v,t)g(x-t v)+I^2(x,v,t)g(x+tv)-2I^1(x,0,0)g(x)}{t}\, d\mu= \nonumber\\
\int_\Omega \Big [ D_1 H (x, Sx, S^2x,...,S^{N-1}x)(v)  + D_1 H (x, Sx, S^2x,...,S^{N-1}x)(-v) \Big] g(x) \, d\mu \leq 0,
\end{eqnarray}
where the last inequality is due to the concavity of $H$ with respect to the first variable.
We shall now prove that the limit of the left hand side of (\ref{fin}) is non-negative as $t \to 0^+$.
It follows from the convexity of $H$ with respect to the last $N-1$ variable together with  $u_i(x) \in \partial_{i+1} H(x, Sx,..., S^{N-1}x)$ that

\begin{eqnarray*}
\int_{\Omega}\frac{I^1(x,v,t)g(x)+I^1(x,-v,t)g(x)-2I^1(x,0,0)g(x)}{t}\, d\mu\geq \\
\frac{1}{t}\int_{\Omega} \sum_{i=1}^{N-1} \langle u_i(x), S^i(x+tv)+S^i(x-tv)-2S(x) \rangle g(x) \, d\mu.
\end{eqnarray*}
The right hand side of the above expression goes to zero, as $t \to 0, $ provided  either $S \in W_{loc}^{1,1}(\Omega)$ or $u_i\in W_{loc}^{1,1}(\Omega)$
for $i=1,2,..., N-1.$
This together with (\ref{fin}) and (\ref{fin1}) imply  that
\begin{equation*}
\int_\Omega \Big[ D_1 H (x, Sx, S^2x,...,S^{N-1}x)(v)  + D_1 H (x, Sx, S^2x,...,S^{N-1}x)(-v)\Big ] g(x)  \, d\mu = 0,
\end{equation*}
from which the desired results follows.\hfill $\square$

\begin{lemma} \label{seq00} For  $v \in \R^d$, define $G_i(v)=\int_\Omega D_iH(x, Sx, ..., S^{N-1}x )(v)\, d\mu.$ Then

 $$\sum_{i=1}^N [G_i(v)+G_i(-v)]\leq 0.$$
\end{lemma}
\textbf{Proof.}
Define $f_i(t,x,v)= H(\sigma^{N+1-i}\big (x+tv, Sx, ..., S^{N-1}x\big))$. Note that $$t \to \frac{f_i(t,x,v)+f_i(t,x,-v)-2f_i(0,x,v)}{t}$$ is monotone and does not change sign.
 It follows from the monotone convergence theorem that
\begin{eqnarray*}\lim_{t \to 0^+}\int_\Omega\frac{f_i(t,x,v)+f_i(t,x,-v)-2f_i(0,x,v)}{t} \, d\mu = \int_\Omega\lim_{t \to 0^+} \frac{f_i(t,x,v)+f_i(t,x,-v)-2f_i(0,x,v)}{t} \, d\mu\\
=
\int_\Omega \big[D_iH(\sigma^{N+1-i}\big (x, Sx, ..., S^{N-1}x\big))(v)+D_iH(\sigma^{N+1-i}\big (x, Sx, ..., S^{N-1}x\big))(v-)\big ] \, d\mu \\
=
\int_\Omega \big[D_iH (x, Sx, ..., S^{N-1}x)(v)+D_iH  (x, Sx, ..., S^{N-1}x)(-v)\big ] \, d\mu=G_i(v)+G_i(-v).
\end{eqnarray*}
Let  $\chi_{\Omega}(t,x)$ be a function that is one when both $x+tv, x-tv \in \Omega$ and zero otherwise. It follows from the dominated convergence theorem that
\begin{eqnarray*}
 G_i(v)+G_i(-v)&=& \int_\Omega\lim_{t \to 0^+} \frac{f_i(t,x,v)+f_i(t,x,-v)-2f_i(0,x,v)}{t} \chi_{\Omega}(t,x) \, d\mu\\
&=&\lim_{t \to 0^+}\int_\Omega\frac{f_i(t,x,v)+f_i(t,x,-v)-2f_i(0,x,v)}{t} \chi_{\Omega}(t,x)  \, d\mu.
\end{eqnarray*}
Let $f(t,x,v)= \sum_{i=1}^N f_i(t,x,v)$. Note that for each $x \in \Omega$ one has $f(t,x,v)= \sum_{i=1}^N f_i(t,x,v) \leq 0$ for $t$ small enough such that $x+tv \in \Omega$.
Similarly $f(t,x,-v) \leq 0$  for $x-tv \in \Omega.$ One also has that $\int_{\Omega}f(0,x,v) \, d\mu=0.$
It follows that
\begin{eqnarray*}
\sum_{i=1}^N  [G_i(v)+G_i(-v)] &=& \int_\Omega\lim_{t \to 0^+} \frac{f(t,x,v)+f(t,x,-v)-2f(0,x,v)}{t} \chi_{\Omega}(t,x)\, d\mu \\
&=&
\lim_{t \to 0^+}\int_\Omega\frac{f(t,x,v)+f(t,x,-v)-2f(0,x,v)}{t} \chi_{\Omega}(t,x) \, d\mu   \\
 &=& \lim_{t \to 0^+}\int_\Omega\frac{f(t,x,v)+f(t,x,v)}{t} \chi_{\Omega}(t,x) \, d\mu  \leq 0.
\end{eqnarray*}
\textbf{Proof of Theorem \ref{dar}.}
Fom Lemma \ref{seq0} and \ref{seq00} we have for  each $v \in \R^d$ and $i=1,2,..., N$
\begin{equation}\label{seq5}
\int_\Omega \Big [ D_i H (x, Sx, S^2x,...,S^{N-1}x)(v)  + D_i H (x, Sx, S^2x,...,S^{N-1}x)(-v) \Big] \, d\mu = 0.
\end{equation}
Since the integrand does not change sign, it has to be zero almost everywhere.
Now choose $\{v_k\}_{k=1}^\infty$ to be a countable dense subset of $\R^d.$
 Set $$A_k=\{  x \in \Omega; D_i H (x, Sx, S^2x,...,S^{N-1}x)(v_k)+D_i H (x, Sx, S^2x,...,S^{N-1}x)(-v_k)=0, \quad 1 \leq i \leq N\}$$
It follows from (\ref{seq5}) that $\Omega\setminus A_k$ is a null set and $\Omega_0=\cap_k A_k$ is  a full measure subset of $\Omega$ such that
$\nabla_i H (x, Sx, S^2x,...,S^{N-1}x)$ exists for all $x \in \Omega_0.$ \hfill $\square$

 \end{document}